\documentstyle[12pt]{article}
\newtheorem{prop}{Proposition}
\newtheorem{def-prop}{Definition-Proposition}
\newtheorem{thm}{Theorem}
\newtheorem{cor}{Corollary}
\newtheorem{conj}{Conjecture}
\newtheorem{lemma}{Lemma}
\newtheorem{defn}{Definition}
\def\endproof{$\Box$}
\def\1{1}
\def\P{{\bf P}}

\def\Z{{\bf Z}}

\def\Q{{\bf Q}}

\def\t{\tau}
\def\Mgbar{\overline{{\cal M}}_g}

\def\liminv{\lim_{\leftarrow V}}
\begin{document}
\title{Equivariant Intersection Theory}
\author{Dan Edidin\thanks{Both authors were partially supported by NSF
postdoctoral
fellowships}\\Department of Mathematics\\University of Missouri\\
Columbia MO 65211\\
William Graham\\ Department of Mathematics\\ University of Chicago\\ Chicago IL
60637}
\date{}
\maketitle
\section{Introduction}
The purpose of this paper is to develop an
equivariant intersection theory for actions
of linear algebraic groups on algebraic schemes. The
theory is based on our construction of equivariant Chow groups. They are
algebraic analogues of equivariant cohomology groups which satisfy
all the functorial properties of ordinary Chow groups. In addition,
they enjoy many of the properties of
equivariant cohomology.
The principal results of this paper are:\\

(1) If a group $G$ acts with finite stabilizers on a
scheme $X$, then rational equivariant Chow groups can be identified
with the rational Chow groups of a quotient.  As a result, we show
that the rational Chow groups of quotients of smooth varieties by
group actions have a canonical ring structure.  This extends and
simplifies previous work of Mumford (\cite{Mu}), Gillet (\cite{Gi})
and Vistoli (\cite{Vi}). In addition the integral Chow groups are an
invariant of the quotient stack $[X/G]$, so we can associate an {\it
integral} Chow ring to smooth quotient stacks.\\

(2) There is a Riemann-Roch isomorphism between a completion
of equivariant $K$-theory of coherent sheaves and a completion
of equivariant Chow groups. This extends the Riemann-Roch theorem
of Baum, Fulton and MacPherson to the equivariant case.\\

(3) There is a localization for torus actions relating the the
equivariant Chow groups of a scheme to the Chow groups of the fixed
locus.  Such a theorem is a hallmark of other equivariant theories
such as cohomology and $K$-theory.  The localization theorems in
equivariant cohomology and $K$-theory imply residue formulas such as
Bott's ( \cite{B-V}, \cite{A-B}, \cite{I-N}),
which can now be proved using intersection theory.

\medskip

Previous work on equivariant intersection theory (\cite{Br},
\cite{Gi}, \cite{Vi}) defined equivariant Chow groups using invariant
cycles on $X$.  The definition we give of equivariant Chow groups, in
contrast, is modeled on Borel's definition of equivariant cohomology.
Borel's insight was to replace the original topological space $X$ by a
homotopic space $X \times EG$, where $EG$ is a contractible space on
which $G$ acts freely.  Since $G$ acts freely, there is a nice
quotient $X_G$ of $X \times EG$ by $G$ and equivariant cohomology is
defined as the cohomology of $X_G$.  To define equivariant Chow groups
one needs an appropriate algebraic replacement for $EG$. This was
supplied by Totaro \cite{To}, who used finite dimensional
representations to approximate the infinite dimensional space $EG$. In
particular if $V$ is a representation of $G$, let $U$ denote an open
set on which $G$ acts freely and has a
a quotient $U \rightarrow U/G$ which is a principal bundle
For any linear algebraic group, the
representation can be chosen so that $V-U$ has arbitrarily large
codimension.  If $X$ is a $G$-scheme then, under mild hypotheses on
$G$ or $X$ (see below), $X \times U$ has a quotient $X \times^G U$ so
that $X \times U\rightarrow X \times^G U$ is a principal $G$-bundle.
The group $A_{\dim V +i - \dim G}(X \times^G U)$ is independent of $V$
as long as the codimension of $V-U$ is sufficiently large. This
defines the $i$-th equivariant Chow group $A_i^G(X)$.

Because $X \times U \rightarrow X \times^G U$ is
a principal $G$-bundle, cycles on $X \times^G U$ exactly
correspond to $G$-invariant cycles on $X \times U$. Since
we only consider cycles of codimension smaller
than the dimension of $X \times (V-U)$, we may in fact
view these as $G$-invariant cycles on $X \times V$.
In other words, instead of considering only $G$-invariant cycles on
$X$ we consider $G$-invariant cycles on $X \times V$ for any
sufficiently big representation of $G$.  By enlarging the class
of cycles we allow, we obtain a theory with many good properties.

By construction, the equivariant Chow groups $A_*^G(X)$ inherit most
of the properties of ordinary Chow groups. In particular if $X$ is
smooth, then there is an intersection product on the equivariant Chow
groups $A_*^G(X)$, no matter how badly $G$ acts on $X$. When $G$ acts
properly and a quotient $X/G$ exists we prove (Theorem \ref{quotient})
that $A_*^G(X)_{\Q}$ = $A_*(X/G)_{\Q}$.  As a result, this proves that
if $G$ acts properly on a smooth variety $X$, then the rational Chow
groups of a quotient $X/G$ have a canonical intersection product
(Corollary \ref{moduli}). This extends the results of Vistoli, who
proved such a theorem when $G$ acts with finite, reduced stabilizers.
This theorem should be useful for doing intersection theory on moduli
spaces of objects which posess infinitesimal automorphisms.  It can
also be used to do intersection theory on toric varieties in arbitrary
characteristic. Furthermore, by avoiding the use of algebraic stacks,
our proof is much simpler.

Another interesting aspect of the theory is that the groups
$A_*^G(X)$ are actually an invariant of the {\em quotient stack}
$[X/G]$ (Proposition \ref{qstacks}).
Thus if $X$ is smooth, then there is an integral intersection
ring associated to the quotient stack $[X/G]$. When $[X/G]$
is Deligne-Mumford (i.e. $G$ acts with finite, reduced stabilizers)
then our ring tensored with $\Q$ agrees with the rings of Gillet
and Vistoli. It would be interesting to compute the torsion
in the equivariant Chow ring in examples of moduli stacks
such as curves of low genus.

Our results on quotient stacks also suggest that there should
be integral Chow rings associated to arbitrary smooth stacks.
Motivated by the equivariant Chow ring, we expect that this ring
would have torsion in arbitrarily high degree. However, we do not
know how to construct such a ring in general.

The connection between equivariant $K$-theory and equivariant Chow groups
is given by our
equivariant Riemann-Roch theorem (Theorem \ref{rockandroll}).
We prove that there is an isomorphism $\t_X: \widehat{K_0^{'G}(X)}_{\Q}
\rightarrow \widehat{A_*^G(X)}_{\Q}$ between the completion of
$K$-theory along the augmentation ideal of the representation ring $R(G)$,
and the completion of the $A_*^G(X)$ along the augmentation
ideal of the equivariant Chow ring of a point.

Along the way we prove a theorem (Theorem \ref{completions})
which shows that the completion of $K_i^{'G}(X)$  along
the augmentation ideal of $R(G)$ is the same as the completion
along the augmentation ideal of $K^0_G(X)$. This result
is related to results of \cite{CEPT} and answers a special
case of a conjecture of \cite{kock}.

In the last part of the paper we prove a localization theorem
for torus actions. If a torus $T$ acts on $X$ with fixed
locus $X^T \subset X$, then $A_*^T(X) \otimes_{R_T} (R_T^+)^{-1} =
A_*(X^T) \otimes (R_T^+)^{-1}R_T$, where $R_T^+$ denotes the multiplicative
set of homogeneous elements of positive degree in $R_T =
Sym(\hat{T})$.

Finally, using the localization theorem, we prove the Bott residue
formula for actions of (split) tori on smooth complete varieties. This
formula has been recently applied in enumerative geometry (cf. \cite{E-S})
so we include an intersection theoretic proof.
Our line of argument follows that of \cite{A-B} using equivariant
intersection theory in place of equivariant cohomology.
(Note that Iversen and Nielsen \cite{I-N} gave an algebraic proof of
this formula  - for smooth projective varieties - using equivariant $K$-theory.
Also, using techniques of algebraic deRahm homology,
H\"ubl and Yekutieli \cite{H-Y}
proved a version - in characteristic 0 - for the action
of an algebraic vector field with isolated fixed points.)

In Section \ref{mixed} we discuss extensions of the theory to
group schemes over a regular base scheme, and
in the appendix we prove a number of technical results about
group actions and quotients in arbitrary characteristic.

{\bf Acknowledgements:} We thank William Fulton,
Rahul Pandharipande and Angelo Vistoli for advice and encouragement.
We also benefitted from discussions with Burt Totaro and Amnon Yekutieli.
Thanks also to Holger Kley for suggesting the inclusion
of the cycle map to equivariant cohomology.
\section{Definitions and basic properties}
\subsection{Conventions and Notation}
Except in Section
\ref{mixed}, all schemes are assumed to be of finite type over a
field of arbitrary characteristic.
A variety is a reduced and irreducible scheme. An algebraic
group is always assumed to be linear.

If an algebraic group $G$ acts on a scheme $X$ then the action is said to
be {\it closed} if the orbits of geometric points are closed in $X$.
It is {\it proper} if the action map $G \times X \rightarrow X \times
X$ is proper. Finally, we say that it is {\it free} if the action map
is a closed embedding. By (\cite[Prop. 0.9]{GIT}) if the action is
free and a geometric quotient scheme $X/G$ exists, then $X$ is
a principal $G$ bundle over $X/G$.

Throughout the paper we will assume that at least one of the following
hypotheses on $X$ or $G$ is satisfied.\\

(1) $X_{red}$ is quasi-projective and the action is linearized
with respect to some projective embedding.\\

(2) $G$ is connected and $X_{red}$ equivariantly embeds as
a closed subscheme in a normal variety.\\

(3) $G$ is special in the sense of \cite{Sem-Chev}; i.e. all
principal $G$-bundles are locally trivial in the Zariski topology.
(Examples of special groups are tori, solvable and unipotent
groups as well as $GL(n)$, $SL(n)$, and $Sp(2n)$. Finite groups
are not special, nor or the orthogonal groups $SO(2n)$ and $SO(2n+1)$.)

\medskip

For simplicity of exposition, we will usually assume that $X$ is
equidimensional.

\subsection{Equivariant Chow groups} Let $X$ be an $n$-dimensional
scheme.
We will denote the $i$-th equivariant Chow
of $X$ group by $A^G_i(X)$. It is defined as follows.

Let $G$ be a $g$-dimensional algebraic group.
Choose an $l$-dimensional representation $V$ of $G$ such that $V$
has an open set $U$ on which $G$ acts freely and whose complement has
codimension more than $n-i$. Assume that a quotient
$U \rightarrow U/G$ (necessarily a principal bundle) exists.
(Such representations exist for any group; see Lemma \ref{q.exist}
of the Appendix.)
The principal bundle $U \rightarrow U/G$ is Totaro's
finite
dimensional approximation of the classifying bundle $EG \rightarrow BG$
(see \cite{To} and \cite{E-G}).
The diagonal
action on $X \times U$ is also free, and since one hypothesis (1)-(3)
holds, there is a quotient
$X_{red} \times U \rightarrow (X_{red} \times U)/G$
which is a principal $G$ bundle\footnote{Without any hypothesis on $X$
or $G$, we only know that the quotient exists as an algebraic space.}
(Prop \ref{inap}).
We will usually denote this quotient
by $(X _{red}\times^G U)$ or $X_G$.

\begin{def-prop} \label{keydef}
Set $A_i^G(X)$ (the $i$-th equivariant Chow group) to be
$A_{i+l-g}(X_G)$, where $A_*$ is the usual Chow group.
This group is independent of the representation
as long as $V- U$ has sufficiently high codimension.
\end{def-prop}

{\bf Remark.} In the sequel, the notation $U \subset V$ will
refer to an open set in a representation on which the action is
free. Because we are working with Chow groups, we
will, when no confusion can arise, abuse notation and act as if all
schemes are reduced.

\medskip

Proof of Definition-Proposition \ref{keydef}.
We will use Bogomolov's double fibration
argument.
Let $V_1$ be another representation of dimension $k$ such
that there is an open $U_1$ with a principal bundle quotient
$U_1 \rightarrow U_1/G$ and whose complement has
codimension at least $n-i$.
Let $G$ act diagonally on $V \oplus V_1$.
Then $V \oplus V_1$ contains an open set
$W$ which has
a principal bundle quotient $W/G$ and contains both
$U \oplus V_1$ and $V \oplus U_1$.
Thus, $A_{i+k+l-g}(X \times^G W) =
A_{i+k+l-g}(X \times^G (U \oplus V_1))$ since
$(X \times^G W)-(X \times^G (U \oplus V_1)$
has dimension smaller than $i+k+l-g$. On the
other hand, the projection $V \oplus V_1 \rightarrow V$
makes $X \times^G (U \oplus V_1)$ a vector
bundle over $X \times^G U$ with fiber $V_1$ and
structure group $G$. Thus, $A_{i+k+l-g}(X \times^G (U \oplus V_1))
= A_{i+l-g}(X \times^G U)$. Likewise,
$A_{i+k+l-g}(X \times^G W) =A_{i+k-g}(X \times^G U_1)$,
as desired.
\endproof

\paragraph{Example} For the classical groups, the representations
and subsets can be constructed explicitly.
In the simplest case, if $G= {\bf G}_m$
then we can take $V$ to an $l$-dimensional
representation with weights one,
$U = V - \{0\}$, and $U/G = \P^{l-1}$. If $G = GL_n$,
take $V$ to be the vector space of $n \times p$
matrices ($p>n$), with $GL_n$ acting by left multiplication, and let
$U$ be the subset of matrices of maximal rank.  Then $U/G$ is the
Grassmannian $Gr(n,p)$.

\paragraph{Remarks}
Now that we have defined equivariant Chow groups,
we will use the notation $X_G$ to mean a mixed quotient $X \times^G U$
for any representation $V$ of $G$. If we
write $A_{i+l-g}(X_G)$ then $V-U$ is assumed to have codimension
more than $n-i$ in $V$. (As above $n=\mbox{dim }X$, $l=\mbox{dim }V$
and $g =\mbox{dim }G$.)

If $Y \subset X$ is an $m$-dimensional $G$-invariant subvariety, then
it has a $G$-equivariant fundamental class $[Y]_G \in A_m^G(X)$.
However, unlike ordinary Chow groups, $A_i^G(X)$
can be non-zero for any $i \leq n$, including negative $i$.
The projection $X \times U \rightarrow U$ induces
a map $X_G \rightarrow U$ with fiber $X$. Restriction
to a fiber gives a map $i^*:A_*^G(X) \rightarrow A_*(X)$
from equivariant Chow groups to ordinary Chow groups. The
map is independent of the choice of fiber because any two
points of $U/G$ are rationally equivalent.
For any $G$-invariant subvariety $Y \subset X$,
$i^*([Y]_G) =[Y]$.

\subsection{Functorial properties}
In this section all maps $f: X \rightarrow Y$ are assumed to
be $G$-equivariant.

If $f: X \rightarrow Y$ is proper, then by descent,
the induced map $f_G: X_G \rightarrow Y_G$
is also proper. Likewise, if $f:X \rightarrow Y$ is flat of
relative dimension $k$ then $f_G:X_G \rightarrow Y_G$ is flat
of dimension $k$.
\begin{defn}
Define proper pushforward $f_*:A_i^G(X) \rightarrow A_i^G(Y)$,
and flat pullback $f^{*}:A_i^G(Y) \rightarrow A_{i-k}^G(X)$
as $f_{G*}:A_{i+l-g}(X_G) \rightarrow A_{i+l-g}(Y_G)$ and
$f_G^*:A_{i+l-g}(Y_G) \rightarrow A_{i+l-g-k}(X_G)$ respectively.
\end{defn}

\medskip

If $f: X \rightarrow Y$ is smooth, then
$f: X_G \rightarrow Y_G$ is also
smooth. Furthermore, if $f: X \rightarrow
Y$ is a regular embedding, then
$f \times id: X \times U \rightarrow Y \times U$
is a regular embedding. In the cartesian diagram
$$\begin{array}{ccc}
X \times U & \rightarrow & Y \times U\\
\downarrow & & \downarrow \\
X_G & \rightarrow & Y_G
\end{array}$$
the vertical arrows are flat and surjective so
by \cite[Prop. IV 3.5]{F-L} the map
$X_G \rightarrow Y_G$ is also a regular embedding.
In particular, if $f: X \rightarrow Y$ is
an l.c.i morphism in
the sense of \cite[Section 6.6]{Fulton},
then $f_G:X_G \rightarrow Y_G$ is also
l.c.i.
\begin{defn} \label{lci}
If $f: X \rightarrow Y$ is l.c.i. of codimension $d$
then define $f^{*}:A_i^G(Y) \rightarrow A_{i-d}^G(X)$
as $f^*_G: A_{i+l-g}(Y_G) \rightarrow A_{i+l-g-d}(X_G)$.
\end{defn}

\begin{prop}
The maps $f_*$ and $f^{*}$ above
are well defined.
\end{prop}
Proof: We will use the double fibration argument.
Let $V_1$ be another representation. Then we have
a cartesian diagram
$$\begin{array}{ccc}
X \times^G (U \oplus V_1) & \rightarrow &Y \times^G (U \oplus V_1)\\
\downarrow & & \downarrow\\
X \times^G U & \rightarrow & Y \times^G U
\end{array}$$
The vertical maps are flat, and their pullbacks
are the isomorphisms which allowed us to define
$A_i^G$. Since flat pullback is compatible with proper
pushfoward, the equivariant pushforward $f_*$ is well defined. Likewise
the flat pullback is compatible with flat and l.c.i
pullback, so $f^{*}$ is also well defined.
\endproof

\medskip

\subsection{Chern classes}
Let $X$ be a scheme with a $G$ action, and let
$E$ be an equivariant vector bundle. For each $i$, $j$
define a map $c_j^G(E):A_i^G(X) \rightarrow A_{i-j}^G(X)$
as follows. Let $V$ be an $l$-dimensional representation such that
$V-U$ has high codimension. By hypothesis, there is
a principal bundle $X \times U \rightarrow X_G$. Thus
by \cite[Prop. 7.1]{GIT} there is a quotient $E_G$ of $E \times U$.

\begin{lemma}
$E_G \rightarrow X_G$ is a vector bundle.
\end{lemma}
Proof.
The bundle $E_G \rightarrow X_G$ is an affine bundle which locally
trivial in the \'etale topology since it becomes locally trivial after the
smooth
base change $X \times U \rightarrow X_G$. Also, the transition
functions are affine since they are affine when pulled back to $X \times U$.
Hence, by descent, $E_G \rightarrow
X_G$
is locally trivial in the Zariski topology and has affine transition
functions; i.e., $E_G$ is a vector bundle over $X_G$.
\endproof

\medskip

Identify $A_i^G(X)$ and $A_{i-j}^G(X)$
with $A_{i+l-g}(X_G)$ and $A_{i-j+l-g}(X_G)$ respectively.
\begin{def-prop}
Define equivariant Chern classes $c_j^G(E):A_i^G(X) \rightarrow A_{i-j}^G(X)$
by $c_j^G(E)\cap \alpha= c_j(E_G) \cap \alpha \in A_{i-j+l-g}(X_G)$.
This definition does not depend on the choice of representation.
\end{def-prop}
Proof: Let $V_1$ be another representation. Then
the pullback of $E \times^G U$ to $X \times^G (U \oplus V_1)$
is isomorphic to the quotient $E \times^G (U \oplus V_1)$.
\endproof

\medskip

Given the above propositions, equivariant Chow groups
satisfy all the formal properties of ordinary Chow groups
(\cite[Chapters 1-6]{Fulton}). In particular, if $X$ is
smooth, there is an intersection product on the
the equivariant Chow groups $A_*^G(X)$ which makes $\oplus A_*^G(X)$
into a graded ring.

\subsection{Operational Chow groups}
Define equivariant operational Chow groups
$A^i_G(X)$ as operations $c(Y \rightarrow X): A_*^G(Y) \rightarrow
A_{*-i}^G(Y)$
for every $G$-map $Y \rightarrow X$.
As for ordinary operational Chow groups (\cite[Chapter 17]{Fulton}),
these operations should be
compatible with the operations on equivariant Chow groups
defined above (pullback
for l.c.i. morphisms, proper pushforward, etc.) From this definition
it is clear that for any $X$, $A^*_G(X)$ has a ring structure.
The ring $A^*_G(X)$  is positively graded, and $A^i_G(X)$ can
be non-zero for any $i \geq 0$.

Note
that by construction, the equivariant Chern classes
defined above are elements of the equivariant operational
Chow ring.

\begin{prop} \label{opsmooth}
If $X$ is smooth of dimension $n$,
then $A^i_G(X) \simeq A_{n-i}^G(X)$.
\end{prop}
Proof. Define a map $A^i_G(X) \rightarrow A_{n-i}^G(X)$
by the formula $c \mapsto c \cap [X]_G$.
Define a map
$A_{n-i}^G(X)\rightarrow A^i_G(X)$, $\alpha \mapsto c_{\alpha}$
as follows.  If $Y \stackrel{f} \rightarrow X$ is a
$G$-map, then since $X$ is smooth, the graph $\gamma_f: Y \rightarrow
Y \times X$ is a $G$-map which is a regular embedding.  If $\beta \in
A^G_*(Y)$ set $c_\alpha \cap \beta= \gamma_f^*(\beta \times
\alpha)$ (note that the cartesian product of equivariant classes is
well defined).

\medskip

Claim (cf. \cite[Proposition 17.3.1]{Fulton}): $\beta \times
(c \cap [X]_G)
 = c \cap (\beta \times  [X_G])$.

\medskip
Given the claim, the formal arguments of
\cite[Proposition 17.4.2]{Fulton}
show that the two maps are inverses.

Proof of Claim: The equivariant class $\beta$ is represented
by a cycle on some $Y_G$, which we can assume to be the fundamental
class of a subvariety $W \subset Y_G$.
Let $\tilde{W}$ be the
inverse image of $W$ in $Y \times U$.
Then $\beta$ pulls back to $[\tilde{W}]_G$
by the equivariant projection map $Y \times U \rightarrow Y$.
By requiring $V - U$ to have sufficiently high codimension, we may
assume that the pullback on Chow groups is an
isomorphism
in the appropriate degrees. Replacing $Y$ by $Y_G$, we
may assume $\beta = [\tilde W]_G$. Since $\tilde W$ is $G$-invariant,
the projection $p: \tilde{W} \times X \rightarrow X$ is equivariant.
Thus,
$$(c \cap [X]_G) \times [\tilde{W}]_G = p^*(c \cap [X]_G)=
c \cap p^*([X]_G) = c \cap ([X_G] \times [\tilde{W}]_G).$$
\endproof

\medskip

Let $V$ be a representation such that
$V- U$ has codimension more than $k$, and set $X_G =
X \times^G U$. In the remainder of the subsection we will
discuss the relation between $A^k_G(X)$ and $A^k(X_G)$ (ordinary
operational Chow groups).

Recall \cite[Definition 18.3]{Fulton} that an envelope
$\pi:\tilde{X} \rightarrow X$ is a proper map such that for any
subvariety $W \subset X$ there is a subvariety $\tilde{W}$ mapping
birationally to $W$ via $\pi$. In the case of group actions, we will
say that $\pi: \tilde{X} \rightarrow X$ is an {\it equivariant} Chow
envelope, if $\pi$ is $G$-equivariant, and if we can take $\tilde{V}$
to be $G$-invariant for $G$-invariant $V$. If there is an open set $X^0
\subset X$ over which $\pi$ is an isomorphism, then we say $\pi:
\tilde{X} \rightarrow X$ is a {\it birational} envelope.

\begin{lemma} If $\pi: \tilde{X} \rightarrow X$ is an
equivariant (birational)  envelope, then
$p: \tilde{X}_G \rightarrow X_G$ is a (birational)
 envelope ($\tilde{X}_G$ and $X_G$ are constructed
with respect to a fixed representation $V$). Furthermore,
if $X^0$ is the open set over which $\pi$ is an isomorphism
(necessarily $G$-invariant), then $p$ is an isomorphism
over $X^0_G = X^0 \times^G U$.
\end{lemma}
Proof: Fulton \cite[Lemma 18.3]{Fulton} proves that
the base extension of an envelope is an envelope.
Thus $\tilde{X} \times U \stackrel{\pi \times id}\rightarrow X \times U$
is an envelope. Since the projection $X \times U \rightarrow X$
is equivariant, this envelope is also equivariant.
If $W \subset X_G$ is a subvariety, let $W'$ be its inverse image
(via the quotient map) in $X \times U$. Let $\tilde{W'}$ be
an invariant subvariety of $\tilde{X} \times U$ mapping
birationally to $W'$. Since $G$ acts freely on $\tilde{X} \times U$
it acts freely on $\tilde{W'}$, and $\tilde{W} = \tilde{W'}/G$
is a subvariety of $\tilde{X}_G$ mapping birationally to $W$.
This shows that $\tilde{X}_G \rightarrow X_G$ is an envelope;
it is clear that the induced map $\tilde{X}_G \rightarrow
\tilde{X}$ is an isomorphism over $X_0^G$. \endproof

\medskip

Suppose $\tilde{X} \stackrel{\pi}\rightarrow X$
is an equivariant envelope which is
an isomorphism over $U$. Let $\{S_i\}$ be the irreducible components
of $S= X -X^0$, and let $E_i = \pi^{-1}(S_i)$. Then $\{S_{i G}\}$
are the irreducible components of $X_G - X^0_G$ and
$E_{i G} = p^{-1}(S_{i G})$.

\begin{thm}
If $X$ has an equivariant non-singular envelope
$\pi: \tilde{X} \rightarrow X$ such that there is an
open $X^0 \subset X$ over which $\pi$ is an isomorphism, then
$A^k_G(X) = A^k(X_G)$.
\end{thm}
Proof: If $\pi: \tilde{X} \rightarrow X$ is an
equivariant non-singular  envelope, then
$p: \tilde{X}_G \rightarrow X_G$
is also an  envelope and $\tilde{X}_G$ is non-singular.
Thus, by \cite[Lemma 1.2]{Kimura}
$p^*:A^*(X_G) \rightarrow A^*(\tilde{X}_G)$ is injective.
The image of $p^*$ is described inductively
in \cite[Theorem 3.1]{Kimura}. A class
$\tilde{c} \in A^*(\tilde{X}_G)$ equals
$p^*c$ if and only if for each
$E_{i G}$ , $\tilde{c}_{| E_{i G}} = p^*c_i$
where $c_i \in A^*(E_i)$.
This description follows from formal properties of operational
Chow groups, and the exact sequence \cite[Theorem 2.3]{Kimura}

$$A^*(X_G) \stackrel{p}\rightarrow A^*(\tilde{X}_G)
\stackrel{p_1^* - p_2^*} \rightarrow A^*(\tilde{X}_G \times_{X_G}
\tilde{X}_G)$$ where $p_1$ and $p_2$ are the two projections
from $\tilde{X}_G \times_{X_G} \tilde{X}_G$.

By Proposition \ref{opsmooth} above, we know that
$A^k_G(\tilde{X}) = A^k(\tilde{X}_G)$.
We will show that $A^k_G(X)$ and $A^k(X_G)$ have the same image
in $A^k(\tilde{X}_G)$.
By Noetherian induction we may assume that
$A^k(S_i) = A^k((S_{i})_G)$. To prove the theorem, it suffices
to show that there is also an exact sequence of equivariant
operational Chow groups
$$0 \rightarrow A^*_G(X) \stackrel{\pi^*}\rightarrow A^*_G(\tilde{X})
\stackrel{p_{1}^* -p_{2}^*}\rightarrow A^*(\tilde{X} \times_X
\tilde{X})$$
This can be checked by working with the action of $A^*_G(X)$
on a fixed Chow group $A_{i}(X_G)$ and arguing as in Kimura's
paper.
\endproof

\begin{cor}
If equivariant resolution of singularities holds (in particular
if the characteristic is 0), and $V-U$ has codimension more than $k$,
then
$A^k_G(X) = A^k(X_G).$
\end{cor}
Proof (c.f. \cite[Remark 3.2]{Kimura}).
We must show the existence of an equivariant envelope
$\pi:\tilde{X} \rightarrow X$. By equivariant
resolution of singularities, there is a resolution
$\pi_1:\tilde{X_1} \rightarrow X$ such
that $\pi_1$ is an isomorphism outside
some invariant subscheme $S \subset X$. By Noetherian
induction, we may assume that we have constructed an
equivariant envelope $\tilde{S} \rightarrow S$. Now
set $\tilde{X} = \tilde{X_1} \cup \tilde{S}$.
\endproof

\subsection{Equivariant higher Chow groups}
In this section assume that $X$ is quasi-projective.
Bloch (\cite{Bl}) defined
higher Chow groups $A^i(X,m)$ as $H_m(Z^i(X,\cdot))$
where $Z^i(X,\cdot)$ is a complex whose $k$-th term
is the group of cycles of codimension $i$ in $X \times \Delta^k$
which intersect the faces properly. Since we prefer
to think in terms of dimension rather than codimension
we will define $A_p(X,m)$ as $H_m(Z_p(X,\cdot))$,
where $Z_p(X,k)$ is
the group of cycles of dimension $p+k$ in $X \times \Delta^k$
intersecting the faces properly. When $X$ is equidimensional
of dimension $n$, then $A_p(X,m) = A^{n-p}(X,m)$.

If $Y \subset X$ is closed, there is a localization long exact sequence.
The advantage of indexing by dimension rather than codimension is that
the sequence exists without assuming that $Y$ is equidimensional.

\begin{lemma}
Let $X$ be equidimensional, and let $Y \subset X$ be closed,
then there is a long exact sequence of higher Chow groups
$$\ldots \rightarrow A_p(Y,k) \rightarrow A_p(X,k) \rightarrow
A_p(X-Y,k) \rightarrow \\
\ldots \rightarrow A_p(Y) \rightarrow A_p(X) \rightarrow A_p(X-Y)
\rightarrow 0$$
(there is no requirement that $Y$ be equidimensional).
\end{lemma}
Proof. This is a simple consequence of
the localization theorem of \cite{Bl}.
By induction it suffices to prove the lemma when $Y$
is the union of two irreducible components $Y_1$, $Y_2$.
In particular, we will prove that the complexes
$Z_p(X-(Y_1 \cup Y_2),\cdot)$ and $\frac{Z_p(X,\cdot)}{Z_p(Y_1 \cup
Y_2,\cdot)}$ are quasi-isomorphic.

By the original localization theorem, $Z_p(X-(Y_1 \cup Y_2),\cdot)
\simeq \frac{Z_p(X-Y_1,\cdot)}{Z_p(Y_2-(Y_1 \cap Y_2),\cdot)}$ and
$Z_p(X-Y_1, \cdot) \simeq \frac{Z_p(X,\cdot)}{Z_p(Y_1)}$.
By induction on dimension, we can assume that the lemma holds for
schemes of smaller dimension, so
$Z_p((Y_2 - (Y_1 \cap Y_2),\cdot) \simeq \frac{Z_p(Y_2,\cdot)}{Z_p(Y_1 \cap
Y_2)}$. Finally
note that
$\frac{Z_p(Y_2),\cdot}{Z_p(Y_1 \cap Y_2)}=\frac{Z_p(Y_1 \cup
Y_2,\cdot)}{Z_p(Y_1,\cdot)}$.
Combining all our quasi-isomorphisms we have
$$Z_p(X-(Y_1 \cup Y_2), \cdot) \simeq \frac{\frac{Z_p(X,\cdot)}
{Z_p(Y_1,\cdot)}}{\frac{Z_p(Y_1 \cup Y_2),\cdot}{Z_p(Y_1,\cdot)}}
\simeq \frac{Z_p(X,\cdot)}{Z_p(Y_1 \cup Y_2,\cdot)}$$
as desired.
\endproof

\medskip

If $X$ is quasi-projective with a $G$-action, we
can define equivariant higher Chow groups
$A_{i}^G(X,m)$ as $A_{i+l-g}(X_G,m)$, where
$X_G$ is formed from an $l$-dimensional representation
$V$ such that $V-U$ has high codimension.
The homotopy lemma for higher Chow groups shows
that $A_{i}^G(X,m)$ is well defined.

Our reason for constructing equivariant higher Chow groups
is to obtain a long exact sequence for a
$G$-invariant subscheme $Y$ of a quasi-projective scheme $X$
with a $G$-action.
\begin{prop} Let $X$ be an equidimensional $G$-scheme, and
let $Y \subset X$ be an invariant subscheme.
There is a long exact sequence of higher equivariant Chow groups
$$\ldots \rightarrow A_p^G(Y,k) \rightarrow A_p^G(X,k) \rightarrow
A_p^G(X-Y,k) \rightarrow \\
\ldots \rightarrow A_p^G(Y) \rightarrow A_p^G(X) \rightarrow A_p^G(X-Y)
\rightarrow 0. $$
\endproof
\end{prop}

\subsection{Cycle Maps}
If $X$ is a complex algebraic variety with the action of
a complex algebraic group, then we can define
equivariant Borel-Moore homology $H_{BM, i}^G(X)$
as $H_{BM,i+2l-2g}(X_G)$ for  $X_G = X \times^G U$.
As for Chow groups, the definition is independent
of the representation, as long as $V -U$ has sufficiently
high codimension,  and we obtain a cycle map
$$cl:A^G_i(X) \rightarrow H_{BM,2i}^G(X)$$
compatible with the usual operations on equivariant
Chow groups (\cite[cf. Chapter 19]{Fulton}).

Note that $H_{BM,i}^G(X)$ is not the same as $H_i(X \times^G EG)$,
where $EG \rightarrow BG$ is the topological classifying bundle.
However, if $X$ is smooth, then $X_G$ is also smooth, and $H_{BM,i}(X_G)$
is dual to $H^{2n-i}(X_G)=H^{2n-i}(X \times^G EG)=H^{2n-i}_G(X)$,
where $n$ is the complex dimension of $X$. In this
case we can interpret the cycle
map as giving a map
$$cl: A^i_G(X) \rightarrow H^{2i}_G(X)$$

If $X$ is compact, and the open sets $U \subset V$ can be chosen so
that $U/G$ is projective, then
Borel-Moore homology of $X_G$ coincides with ordinary
homology, so $H^G_{BM*}(X)$ can be calculated with a compact model.
However
In general, however, $U/G$ is only quasiprojective.
If $G$ is finite, then $U/G$ is
never projective.  If $G$ is a torus, then $U/G$ can be taken to be a
product of projective spaces.  If $G = GL_n$, then $U/G$ can be taken
to be a Grassmannian (see the example in Section 2.1)

If $G$ is semisimple, then $U/G$ cannot be chosen
projective, for then the hyperplane class would be a nontorsion
element in $A^1_G$, but by \cite{E-G}
$A^*_G \otimes \Q \cong S(\hat{T})^W \otimes
\Q$, which has no elements of degree 1.  Nevertheless for
semisimple (or reductive) groups we can obtain a cycle map
$$cl: A_*^G(X)_{\Q} \rightarrow H_{BM*}^T(X;\Q)^W$$ by identifying
$A_*^G(X) \otimes \Q$ with $A_*^T(X)^W \otimes \Q$ and
$H_{BM*}^G(X;\Q)$ with $H_{BM*}^T(X;\Q)^W$; if $X$ is compact then
the last group can be calculated with a compact model.

\section{Intersection theory on quotients}
One of the uses of equivariant intersection theory is
to study intersection theory on quotient stacks
and their moduli. In particular, we show below that the rational
Chow groups of moduli spaces
which are group quotients of a smooth variety
have an intersection product -- even when there are
infinitesimal automorphisms.

\subsection{Chow groups of quotients}
Let $G$ act on a scheme $X$, and assume
that a geometric quotient\footnote
{In characteristic $p$, the definition of geometric quotient used
here is slightly weaker
than the one given in \cite{GIT}. See the appendix.}
$X \rightarrow X/G$ exists.

\begin{prop} \label{p.quotient}
If $G$ acts freely, then $A_*^G(X,m) = A_*(X/G,m)$ (the isomorphism
of higher Chow groups requires $X$ to be quasi-projective).
\end{prop}
Proof. If the action is free, then $(V \times X)/G$ is a vector bundle
over $X/G$. Thus $X_G$ is an open set in this bundle
with arbitrarily high codimension, and the proposition follows from homotopy
properties of (higher) Chow groups. \endproof

\medskip

\begin{thm} \label{quotient}
If $G$ acts properly
on a quasi-projective variety
$X$, so that $X/G$ is quasi-projective,
then\\

(1)$A_*^G(X,m) \otimes \Q = A_*(X/G,m) \otimes \Q$ for all $m \geq 0$.\\

(2) There is an isomorphism of operational Chow rings
$$p^*:A^*(X/G)_{\Q} \stackrel{\simeq} \rightarrow A^*_G(X)_{\Q}.$$
\end{thm}

{\bf Remarks.} (1) If the action is proper,
then the stabilizers are complete.
Since $G$ is affine they must in fact be finite.
We will
sometimes mention that the stabilizers are finite for emphasis.

(2) When $X$ is the set of stable points for some
linearized action of $G$ and $X/G$ is quasi-projective, then the
action is proper.

(3) The condition that $X/G$ is quasi-projective is only
required because the localization long exact sequence for higher Chow
groups has only been proved for quasi-projective schemes.

\medskip

In practice, many interesting varieties arise as quotients
of a smooth variety by a connected
algebraic group
which acts with finite stabilizers.
Examples include simplicial toric varieties
and various moduli spaces such as curves, vector bundles, stable maps, etc.
There is a long history
of work on the problem of constructing an intersection product
on the rational Chow group of quotients of smooth varieties.
In characteristic 0,
Mumford (\cite{Mu}) proved the existence
of an intersection product on the rational Chow groups
of $\Mgbar$, the moduli space of stable curves.
Gillet (\cite{Gi}) and Vistoli (\cite {Vi}) constructed intersection
products on quotients in arbitrary characteristic -- provided
that the stabilizers of geometric points are reduced. In characteristic 0,
Gillet (\cite[Thm 9.3]{Gi}) showed that his product on $\Mgbar$ agreed
with Mumford's, and in \cite[Lemma 1.1]{Edidin} it was shown
that Vistoli's product also agreed with Mumford's.

As a corollary to Theorem \ref{quotient} we
obtain a simple proof of the existence of intersection products on the
rational Chow groups of quotients for a group acting with
finite but possibly non-reduced stabilizers. Furthermore, when the
stabilizers are reduced,
our
product agrees with Gillet's and Vistoli's (Proposition \ref{triprod}).
In particular, this
answers \cite[Conjecture 6.6]{Vi}
affirmatively for moduli spaces of quotient stacks.

\begin{cor} \label{moduli}
Let $Y$ be a quasi-projective variety
which is isomorphic to a geometric quotient
$X/G$, where $X$ is smooth and $G$ acts properly
(hence with finite stabilizers) on $X$.  Then the rational Chow groups
$A_*(Y)_{\Q}$ have an intersection product.  This product is
independent of the presentation of $Y$ as a quotient.
\end{cor}

Proof of Corollary \ref{moduli}. Since $X$ is smooth, the equivariant
Chow groups $A_*^G(X)$ have an intersection product induced by the
isomorphism $A^*_G(X) \rightarrow A_*^G(X)$ (with $\Z$ coefficients).
By Theorem \ref{quotient} $A_*^G(X)_{\Q} = A_*(Y)_{\Q}$ so the groups
$A_*(Y)_{\Q}$ inherit a ring structure.  Since the intersection
product on $A_*(Y)$ is induced by the multiplication in $A^*(Y)$ it
depends only on $Y$.
\endproof

\medskip

Proof of part (1) of Theorem \ref{quotient}.
For simplicity of exposition we give the proof
assuming that the group $G$ is connected of
dimension $g$. (This way we can assume that the set-theoretic
inverse image in $X$  of a subvariety of $X/G$ is a single variety rather than
a possible disjoint union of varieties.)
All coefficients -- including those of cycle groups -- are assumed to
be rational.
If $G$ acts properly on $X$, then $G$ acts
properly on $X \times \Delta^m$
by acting trivially on the second factor. In this case,
the boundary map of the higher
Chow group complex preserves invariant cycles, so there is
a subcomplex of invariant cycles $Z_*(X,\cdot)^G$.
Set $$A_*([X/G],m) = H_m(Z_*(X, \cdot)^G,\partial).$$

Now if $X \rightarrow X/G$ is a geometric quotient, then so
is $X \times \Delta^m \stackrel{\pi} \rightarrow X \times \Delta^m$.
Define a map
$\pi^*: Z_k(X,m) \otimes \Q \rightarrow Z_{k+g}(X,m)^G
\otimes \Q$ for all
$m$ as follows. Let $F \subset X/G  \times \Delta^m$ be a $k+m$-dimensional
subvariety
intersecting the faces properly,
then $H = (\pi^{-1}F)_{red}$  is a $G$-invariant
$(k+m+g)$-dimensional subvariety of $X \times \Delta^m$ which intersects the
faces properly. Thus, $[H] \in Z_{k+g}(X,m)^G $.
Let $e_H$ be the order of the stabilizer
at a general point of $H$, and let $i_H$ be the degree of the purely
inseparable
extension $K(F) \subset K(H)^G$. Set $\pi^*[F] = \alpha_H[
H] \in Z_{k+g}^G(X,m)$, where $\alpha_H = \frac{e_H}{i_H}$.
Since $G$-invariant subvarieties of $X \times \Delta^m$
exactly correspond to subvarieties of $X/G \times \Delta^m$,
$\pi^*$ is an isomorphism of cycles for all $m$.

\begin{prop} \label{proper}
Let
$$\begin{array}{ccc}
Z & \stackrel{g} \rightarrow & X \\
\small{p} \downarrow & & \small{\pi} \downarrow\\
Q & \stackrel{f} \rightarrow & Y
\end{array}$$
be a commutative diagram of quotients with $f$ and $g$ proper.
Then $p^* \circ g_* = f_* \circ \pi^*$ as maps $Z_*(Q) \rightarrow
Z_*(Z,m)^G$.
\end{prop}
Proof of Proposition \ref{proper}. The proposition is an
immediate consequence of the following lemma.

\begin{lemma} Suppose $G$ acts properly (hence with finite
stabilizers) on varieties $Z$ and $X$. Let
$$\begin{array}{ccc}
Z & \stackrel{g} \rightarrow &X \\
\small{p} \downarrow & & \small{\pi} \downarrow\\
Q & \stackrel{f} \rightarrow & Y
\end{array}$$
be a commutative diagram of geometric quotients with
$f$ and $g$ finite and surjective.  Then
$$\frac{e_Z}{i_Z} [K(Q):K(Y)] = \frac{e_X}{i_X} [K(Z):K(Z)].$$
\end{lemma}
Proof. Since we are checking degrees, we may replace $Y$ and $Q$
by $K(Y)$ and $K(Q)$, and $X$ and $Z$ by their generic fibers over
$Y$ and $Q$ respectively.
Then we have a commutative diagram of varieties.
$$\begin{array}{ccc}
Z &  \rightarrow & X \\
\downarrow & &  \downarrow\\
\mbox{spec}(K(Z)^G) & \rightarrow & \mbox{spec}(K(X)^G)\\
\downarrow & &  \downarrow\\
\mbox{spec}(K(Q)) & \rightarrow & \mbox{spec}(K(Y))
\end{array}$$
Since $i_Z := [K(Z)^G:K(Q)]$ and $i_X :=[K(X)^G:K(Y)]$, it suffices
to prove
$$e_Z[K(Z)^G:K(X)^G] = e_X[K(Z):K(X)].$$

By \cite[Paragraph 6.5]{Borel} the
extensions $K(Z)^G \subset K(Z)$ and $K(X)^G \subset
K(X)$ are separable (transcendental).
Thus, after  finite separable base extensions,
we may assume that there are sections $s: \mbox{spec}(K(X)^G) \rightarrow
U$ and $t: \mbox{spec}(K(Z)^G) \rightarrow
W$. In this case the stabilizer group schemes over $W$ and $U$
are
isomorphic to $K(Z)^G \times G$ and $K(X)^G \times G$ respectively
(\cite[Proof of Prop 0.7]{GIT}). Thus, $e_Z = [K(Z)^G \times K(G):
K(Z)]$ and $e_X = [K(X)^G \times K(G): K(X)]$. The lemma follows.
\endproof

\begin{prop} \label{squiggy}
The map $\pi^*$ commutes with the boundary operator of the higher Chow groups.
In particular, there is an induced isomorphism of Chow groups
$$A_k(X/G,m) \simeq A_{k+g}([X/G],m).$$
\end{prop}
Proof of Proposition \ref{squiggy}.
If
$$\begin{array}{ccc}
Z & \stackrel{g} \rightarrow & X \\
\small{p} \downarrow & & \small{\pi} \downarrow\\
Q & \stackrel{f} \rightarrow & X/G
\end{array}$$
is a commutative diagram of quotients with $f$ and $g$ finite and
surjective, then $f_*$ and $g_*$  are surjective as
maps of cycles. Thus, by Proposition \ref{proper} it suffices
to prove $p^*:Z_*(Q) \rightarrow Z_*(X)^G$ commutes with
$\partial$.

By Proposition \ref{whizzbang} of the appendix, there is
a commutative diagram of quotients such that $p:Z \rightarrow Q$
is a principal bundle. Since $p$ is flat, $p^*$ commutes
with $\partial$ and Proposition \ref{squiggy} follows.
\endproof

\medskip

Suppose $T \subset X$ is a $G$-invariant subvariety. Let $S \subset X/G$
be its image under the quotient map. Set $U = X -T$ and $V=X/G -U$.
Then we have two commuative diagrams of geometric quotients

$\begin{array}{ccc} T & \stackrel{i} \rightarrow & X\\
\small{\pi} \downarrow & & \small{\pi} \downarrow \\
S & \stackrel{i} \rightarrow & X/G
\end{array}$ \hspace{2.0in}
$\begin{array}{ccc} U & \stackrel{j} \rightarrow & X\\
\small{\pi} \downarrow & & \small{\pi} \downarrow \\
V & \stackrel{j} \rightarrow & X/G
\end{array}$

\begin{lemma} \label{iggypop}
Let $\alpha \in Z_k(X/G,m)$ and $\beta \in Z_k(S,m)$.

(1) $\pi^*j^* \alpha = j^* \pi^* \alpha$ in $Z_{k+g}^G(U,m)$.

(2) $\pi^*i_* \beta = i_*\pi^*\beta$ in $Z_{k+g}^G(X,m)$.
\end{lemma}

Proof of Lemma \ref{iggypop}. If $\alpha = [F]$ and $H= \pi^{-1}(F)_{red}$,
then $\pi^*j^*\alpha$ and $j^*\pi^*\alpha$ are both multiples of
$[H \cap U]$. Since $e_{[H \cap U]} = e_{[H]}$, and $i_{[H \cap U]} = i_{[H]}$,
the multiplicities are the same. This proves (1).

Part (2) was proved in Proposition \ref{proper}. \endproof

As a consequence of Proposition \ref{squiggy} and Lemma \ref{iggypop},
we obtain
the following proposition.
\begin{prop} \label{snodgrass}
Let $T \subset X$ be an invariant subvariety. If $S$, $U$, and $V$
are as above, then there is a commutative diagram of isomorphisms
$$\begin{array}{ccccccc} \ldots \rightarrow & A_*([T/G],m) & \rightarrow &
A_*([X/G],m) & \rightarrow & A_*([U/G],m) & \rightarrow \ldots \\
& \simeq \uparrow & & \simeq \uparrow & & \simeq \uparrow & \\
\ldots \rightarrow & A_*(S,m) & \rightarrow &
A_*(X/G,m) & \rightarrow & A_*(V,m) & \rightarrow \ldots
\end{array}$$
\endproof
\end{prop}

Next, note that there is a map
$$\alpha:A_*([X/G],m) \rightarrow A_*^G(X,m)$$
defined by the formula
$$[F] \in Z^G_*(X,m) \mapsto [F]_G$$ which commutes with equivariant
proper pushforward and equivariant flat pullback.
\begin{prop} \label{warhol}
If $G$ acts properly
on $X$, and a quasi-projective
geometric quotient $X \rightarrow X/G$ exists,
then $\alpha$ is an isomorphism.
\end{prop}

Proof of Proposition \ref{warhol}.
By the Nullstellensatz there is a point of
$X$ which is finite over the generic point of $X/G$.
Thus, by generic flatness, there is a locally closed
subvariety $Z \subset X$, and an open set $W \subset X/G$ such
that the projection $Z \rightarrow W$ is finite and flat. Let $U =
\pi^{-1}(W)$.
Since $G$ acts properly, the map
$G \times Z \rightarrow U$ is finite. Shrinking $Z$ (and thus $U$)
we may assume that $G \times Z \rightarrow U$ is also flat.
By Noetherian induction and the localization long exact sequence
(which exists by Proposition \ref{snodgrass})
it suffices to prove that $\alpha:A_*^G([U/G],m) \rightarrow  A_*^G(U,m)$
is an isomorphism.

Taking Chow groups, we obtain a commutative diagram where all maps commute.
$$\begin{array}{ccc}
A_*^G(G \times Z,m) & \stackrel{\leftarrow} \rightarrow & A_*^G(U,m)\\
\small{\alpha} \downarrow & & \small{\alpha} \downarrow \\
A_*([G \times Z/G],m) & \stackrel{\leftarrow} \rightarrow & A_*([U/G])
\end{array}$$
(The right horizontal arrows are proper pushforward divided by the degree,
and the left horizontal arrows are flat pullback.)
Chasing the diagram shows that if the left vertical arrow is an isomorphism
then so is the right vertical arrow. Since $G$ acts freely,
$\alpha: A_*([G \times Z/G],m) \rightarrow A_*(G \times Z, m)^G$ is
an isomorphism by Proposition \ref{p.quotient}.
\endproof.

We have now proved part (1) of Theorem \ref{quotient}.

\medskip

Proof of part (2) of Theorem \ref{quotient} (cf. \cite[Proposition 6.1]{Vi}).
Suppose $c \in A^*(X/G)_{\Q}$,
$Z \rightarrow X$ is a $G$-equivariant morphism, and $\alpha \in A_*^G(Z)$.
For any representation $V$, there are maps
$Z_G \rightarrow X_G \rightarrow X/G$. If $V$ is chosen so that
$\alpha$ is represented by a class $\alpha_V \in A_*(Z_G)$ we can define
$$p^*c \cap \alpha = c \cap \alpha_V. $$
As usual, this definition is independent of the representation,
so $p^*c \cap \alpha \in A_*^G(Z)$.

(1) $p^*$ is injective.

Proof of (1). Suppose $p^* \cap \alpha =0$ for all $G$-maps $Z \rightarrow
X$ and all $\alpha \in A_*^G(Z)$. By base change, it suffices to show
$c \cap x = 0$ for all $x \in A_*(X/G)_{\Q}$. By Proposition \ref{whizzbang}
there is a finite map $Y \rightarrow X/G$, and  a principal bundle
$Z \rightarrow Y$ together with a finite $G$-map $Z \rightarrow X$.
Thus we obtain a commutative diagram
$$\begin{array}{ccc}  Z_G & \stackrel{g} \rightarrow & X_G\\
\small{q}\downarrow & & \small{p} \downarrow \\
Y & \stackrel{f} \rightarrow & X/G
\end{array}$$
where the horizontal maps are proper and surjective.
Choose $y \in A_*(Y)$ so that $f_*(y) = x$.
Since $q$ is flat
$$0= c \cap q^*y  =q^*(c \cap y).$$
Since $q^*$ is an isomorphism in the appropriate degrees,
$c \cap y = 0$. Thus
$$0 = f_*(c \cap y) = c \cap f_*y = c \cap x$$
as desired.

(2) $p^*$ is surjective.

Proof of (2). Suppose $d \in A^*_G(X)$. Define $c \in A^*(X/G)$
as follows: If $Y \rightarrow X/G$ and $y \in A_*(Y)$, set
$c \cap y = d \cap \pi^*y$ where $\pi: X\times_{X/G} Y \rightarrow
Y$ is the quotient map, and $\pi^*y \in A_*([X\times_{X/G} Y/G])_{\Q}
= A_*^G(X \times_{X/G} Y)_{\Q}$ is defined as above.
Then $p^*c = d \in A^*_G(X)$.
\endproof.

\subsection{Intersection products on moduli}
Equivariant intersection theory gives a nice way of working with
cycles on a singular moduli space ${\cal M}$ which is a quotient $X/G$ of
a smooth variety by a group acting with finite stabilizers.
Given a subvariety $W \subset {\cal M}$ and
a family $Y \stackrel{p}\rightarrow B$ of schemes
parametrized by ${\cal M}$ there is a map $B \stackrel{f} \rightarrow
{\cal M}$. We wish to define
a class $f^*([W]) \in A_*B$ corresponding to how the image
of $B$ intersects $W$. This can be done (after tensoring with $\Q$) using
equivariant theory.

By Theorem \ref{quotient}, there is an isomorphism
$A_*({\cal M})_{\Q}= A_*^G(X)$ which takes
$[W]$ to the equivariant class $w= \frac{e_W}{i_W}[f^{-1}W]_G$.
Let $B_G \rightarrow B$ be the principal $G$-bundle $B \times_{[X/G]} X$,
(The fiber product is a scheme, although the product is taken over
the quotient stack $[X/G]$.
Typically, $B_G$ is the structure bundle of
a projective bundle $\P(p_*L)$ for a relatively very ample line bundle
$L$ on $Y$).
Since $X$ is smooth, there is an equivariant pullback
$f^*G: A_*^G(X) \rightarrow A_*^G(B_G)$ of the induced map $B_G \stackrel{f_G}
\rightarrow X$, so we
can define a class $f_G^*(w) \in A_*^G(B_G)$. Identifying
$A_*^G(B)$ with $A_*(B)$ we obtain our class $f^*(W)$.

\paragraph{Example (Moduli of stable curves)} Let $\Mgbar$ be the moduli space
of curves. Let $\Delta_{i} \subset
\Mgbar$ be the Weil divisor corresponding
to (stable) nodal curves which are formed by identifying a curve of genus $i$
to
a curve of genus $g-i$ at a point. Given a family of curves
$Y \stackrel{p} \rightarrow B$ there is a map
$B \rightarrow \Mgbar$. We wish to define a cycle $\delta_B$
corresponding to the intersection of the image of $B$ with $\Delta_i$.
Such a class can be defined using Vistoli's intersection theory on the
Deligne-Mumford stack $F_{\Mgbar}$ since there is a Gysin pullback
$A_*(B) \otimes \Q \stackrel{i^*}\rightarrow A_*B \otimes \Q$
corresponding to the inclusion $\delta_i \stackrel{i} \hookrightarrow
F_{\Mgbar}$.
Let $\delta_i$
be the inverse image of $\Delta_i$ in $F_{\Mgbar}$.
If the image of $B$ is completely contained
in $\Delta_i$, then, to calculate $\delta_B$ we need to use an
excess intersection
formula
$$\delta_B = i^*(\delta_i) = c_1({\cal O}_{\delta_i}(\delta_i)) \cap [B].$$
Similarly, if $\Theta \subset \Mgbar$ is a subvariety of codimension
$d$ corresponding
a smooth substack $\theta \subset F_{\Mgbar}$, then we
would like to assert that
$$\theta_B = c_d(N_{\theta}F_{\Mgbar}) \cap [B]$$
where $N_{\theta}F_{\Mgbar}$ is the normal bundle to
$\theta$ in $F_{\Mgbar}$.
Unfortunately, such formulas were not fully developed in \cite{Vi},
so their use can not be completely justified.

Consider the equivariant situation.
The moduli space $\Mgbar$ is a geometric quotient $H_g/G$
where $H_g$ is the (smooth) Hilbert scheme of pluricanonically
embedded stable curves and $G=PGL(N)$ for some $N$.
Given an abstract family $Y \rightarrow B$
let $Y_H \rightarrow B_H$ be the corresponding family of pluricanonically
embedded curves. Let $\Delta_i^H \stackrel{i_H} \hookrightarrow
H_g$ be the corresponding $G$ invariant divisor in $A_*^G(H_g)$.
Since $H_g$ is smooth, we obtain an equivariant line bundle
${\cal O}(\Delta_i^H)$.
Let $B_H \rightarrow H_g$ be the corresponding
equivariant map. Then by equivariant  excess intersection (which follows
from ordinary excess intersection on schemes $(H_g \times U)/G$)
$$i_H^*(\Delta_i^H) = c_1({\cal O}_{\Delta_i^H}(\Delta_i^H)) \cap [B].$$
The corresponding formula in $A_*B \otimes \Q$ follows from
the identification of $A_*(B)$ with $A_*^G(B_H)$.

Similar formulas
(involving Chern classes of the equivariant normal bundle) hold
for cycles of higher codimension.
In \cite[Section 3]{Edidin} explicit excess intersection formulas were given
in codimension 1 and 2, for various nodal loci. The approach
there is similar to the discussion above, although
equivariant Chow groups were not used. Instead, a result of Vistoli (proved
in characteristic 0, although true in arbitrary characteristic)
was used to identify $A^1(B)$ (codimension-one cycles) with
$A^1(B_H)$ because $B_H \rightarrow B$ is a principal $\P GL(N)$ bundle.
The statements of \cite[Section 3]{Edidin} can proved in
in arbitrary characteristic using the methods outlined above.

\subsection{Chow groups of quotient stacks}
If $G$ acts on $X$ we let $[X/G]$ denote the quotient stack.  This is
a stack in the sense of Artin, and exists without any assumptions on
the $G$-action.  By the next proposition, the equivariant Chow groups
do not depend on the presentation as a quotient, so they are an invariant
of the stack.

\begin{prop} \label{qstacks}
Suppose that $[X/G] \simeq [Y/H]$ as quotient stacks. Then
$A_i^G(X) \simeq A_i^H(Y)$ for all $i$.
\end{prop}
Proof: Suppose $\mbox{dim } G =g$ and $\mbox{dim } H = h$. Let
$V_1$ be an $l$-dimensional representation of $G$, and $V_2$ an
$M$ dimensional representation of $H$. Let $X_G = X \times^G U_1$
and $Y_H = X \times^H U_2$, where $U_1$ (resp. $U_2$) is an open set
on which $G$ (resp. $H$) acts freely.
Since the diagonal of a quotient stack is representable,
the fiber product $Z=X_G \times_{[X/G]} Y_H$ is a scheme. This
scheme is a bundle over $X_G$ and $Y_H$ with fiber $U_2$ and $U_1$
respectively.
Thus, $A_{i+l-g}(X_G) = A_{i+l+m-g}(Z) = A_{i+m-h}(Y_H)$
and the proposition follows.
\endproof

\medskip

\noindent
{\bf Remark.} Proposition \ref{qstacks} suggests that there should be a
notion of
Chow groups of an arbitrary algebraic stack which can have non-zero
torsion in arbitrarily high degree. This situation would be analogous to
the cohomology of quasi-coherent sheaves on the \'etale (or flat) site
(cf. \cite[p. 101]{D-M}).

\medskip

If $G$ acts properly with finite, reduced stabilizers, then $[X/G]$ is
a separated
Deligne-Mumford stack. The rational Chow groups $A_*([X/G]) \otimes
\Q$ were first defined by Gillet \cite{Gi} and coincide with the groups
$A_*([X/G]) \otimes \Q$ defined above.
More generally, if $G$ acts with finite stabilizers which are not reduced, then
then Gillet's definition
can be extended and we can define the ``naive'' Chow groups
$A_k([X/G])_{\Q}$ as the group generated by $k$-dimensional integral
substacks modulo rational equivalences.  With this definition we
expect that $A_*^G(X)_{\Q} = A_*([X/G])_{\Q}$. To prove such an
isomorphism in general requires that the naive Chow groups of the stack satisfy
the homotopy property (i.e. if $F \rightarrow G$ is a vector bundle in
the category of stacks, then $A_*(F)_{\Q} = A_*(G)_{\Q})$.
However, if a quasi-projective
quotient exists, then Proposition \ref{warhol} can be restated in
the language of stacks as
\begin{prop}
Let $G$ be a $g$-dimensional group which acts
properly on a scheme $X$
(so the quotient $[X/G]$ is a separated Artin stack). Assume that
a quasi-projective moduli scheme $X/G$ exists for $[X/G]$.
Then $A_i^G(X) \otimes \Q = A_{i-g}([X/G]) \otimes \Q$.
\endproof
\end{prop}

{\bf Remarks.} (1) Although $A_*^G(X) \otimes \Q = A_*([X/G]) \otimes \Q$,
the integral Chow groups may have non-zero torsion for all $i < \mbox{dim }
X$. It would be interesting to compute this torsion in examples
such as moduli spaces of curves of low genus.

(2) In general, a separated quotient stack should always
have an algebraic space as coarse moduli space. Thus, an extension
of the present theory to algebraic spaces would eliminate the need
for any assumptions in the proposition.

\medskip

With the identification of $A_*^G(X) \otimes \Q$ and $A_*([X/G]) \otimes \Q$
there are three intersection products on the rational Chow groups
of a smooth Deligne-Mumford quotient stack with a moduli
space -- the equivariant product, Vistoli's product
defined using a via a gysin pullback for regular embeddings of stacks,
and Gillet's product defined using the product in higher $K-$theory.
The next proposition shows that they are identical.

\begin{prop} \label{triprod}
If $X$ is smooth and $[X/G]$ is a separated
Deligne-Mumford stack
(so $G$ acts with finite, reduced stabilizers)
with a quasi-projective moduli space $X/G$, then
the intersection products on $A_*([X/G])_{\Q}$ defined by Vistoli and
Gillet are the same as the equivariant product on $A_*^G(X)_{\Q}$.
\end{prop}

Proof. If $V$ is an $l$-dimensional representation, then all three products
agree on the smooth quotient scheme (\cite{Vi}, \cite{Grayson})
$(X \times U)/G$. Since
the flat pullback of stacks $f:A^*([X/G])_{\Q} \rightarrow
A^*((X \times U)/G)_{\Q}$ commutes with all 3 products, and is an isomorphism
to arbitrarily high codimension, the proposition follows.
\endproof

\medskip

\section{Equivariant Riemann-Roch}
In this section we construct an equivariant Todd class map and prove
an equivariant Riemann-Roch theorem for $G$-schemes.
The theorem involves
completions of equivariant $K$-groups and Chow groups
because the groups $A_*^G(X)$ (resp. $A^*_G(X)$) can have terms of
arbitrarily large negative (resp. positive) degree. Thus,
the Todd class and Chern character map to completions of
these groups.  The map $\t^G_X: K^{'G}_0(X)
\rightarrow \widehat{A_*^G(X)_{\Q}}$ factors through a completion
of $K_{0}^{'G}(X)$ and we obtain an isomorphism
$\t^G_X: \widehat{K^{'G}_0(X)} \rightarrow \widehat{A_*^G(X)_{\Q}}$.

This section has two parts.  In the first, we define
$\widehat{A_*^G(X,i)}$ and $\widehat{K^{'G}_i(X)}$ as completions of
$A_*^G(X,i)$ and $K^{'G}_i(X)$ along certain ideals.  We then prove an
analogue of a theorem of Atiyah and Segal which gives a more geometric
description of these completions.  In the second part we construct
the Todd class map $\tau_X^G:K_0^{'G}(X) \rightarrow \widehat{A_*^G(X)_{\Q}}$
and show that it induces an isomorphism
$\tau_X^G: \widehat{K^{'G}_i(X)} \stackrel{\simeq} \rightarrow
\widehat{A_*^G(X)_{\Q}}$.
The construction is an easy consequence of the
nonequivariant Riemann-Roch theorem and our geometric
description of the completions.  Finally, we discuss a
conjecture of Vistoli.

\subsection{Completions of equivariant K-groups and Chow groups} Let
$R(G)$ denote the representation ring of $G$.
Let $K_0^G(X)$ denote the Grothendieck group of
$G$-equivariant vector bundles on $X$, and let $K_i^{'G}(X)$ denote
the $i$-th higher $K$-group of the category of $G$-equivariant
coherent sheaves (\cite{Thomason}). As in the
non-equivariant case, $K_0^G(X)$ is a ring under tensor product, and
$K_i^{'G}(X)$ is a module for that ring. Also, $K_0^G(X)$ and
$K_i^{'G}(X)$ are $R(G)$ modules via the isomorphism $R(G) \simeq
K_0^G(pt)=K_0^G$.

Let $P \subset K_0^G = R(G)$ denote the ideal of virtual
representations of dimension 0, and let
$\widehat{K_i'^G(X)}$ be the completion  of $K_i'^G(X)$ along $P$.
Let $Q = A^+_G
\subset A^*_G$ be the augmentation ideal, and let
$\widehat{A_*^G(X,i)}$ be the completion of $A_*^G(X,i)$
along $Q$.

Let $\tilde{Q} = A^+_G(X) \subset A^*_G(X)$ denote the augmentation
ideal, then $Q A^*_G(X) \subset \tilde{Q}$. Let
$\widetilde{A_*^G(X)}$ denote the completion of $A^*_G(X)$ along
$\tilde{Q}$.  Likewise, let  $\tilde{P} \subset K_0^G(X)$ denote the ideal of
virtual bundles of rank 0 (the kernel of the rank map), then $P
K_0^G(X) \subset \tilde{P}$.  Let  $\widetilde{K_i'^G(X)}$ denote the
completion of $K_i'^G(X)$ along $\tilde{P}$.

We will show below that there are isomorphisms
$$
\widetilde{K_i'^G(X)} \cong \widehat{K_i'^G(X)}
$$
$$
\widetilde{A_*^G(X,i)} \cong \widehat{A_*^G(X,i)}.
$$
To show this, we will compare these completions with more geometrically
defined ones.

Partially order the set ${\cal V}$ of representations of $G$ by the rule
$W < V$ if $W$ is a summand in $V$. For each representation,
let $V^f$ be the open set of points whose orbits
are closed in $V$ and which have trivial stabilizer.
The collections
$\{K_i^{'G}(X \times V^f)\} _{ V \in {\cal V} }$ and
$\{A_*^G((X \times V^f),i)\}_{ V \in {\cal V} }$
are inverse systems since the inclusion $V^f \oplus W \hookrightarrow
(V\oplus W)^f)$ induces restriction maps
$$
K_i^{'G}(X \times (V\oplus W)^f) \rightarrow K_i^{'G}(X \times V^f \oplus W)
\simeq K_i^{'G}(X \times V^f)
$$
$$
A_*^{G}((X \times (V\oplus W)^f),i) \rightarrow
A_*^{G}((X \times V^f \oplus W),i)
\simeq A_*^{G}((X \times V^f),i).
$$
By identifying $K_i^{'G}(X \times V)$ with $K_i^{'G}(X)$ and
$A_*^{G}((X \times V),i)$ with $A_*^G(X,i)$ we obtain
restriction maps
$$r'_V:K_i^{'G}(X) \rightarrow K_i^{'G}(X \times V^f)$$
$$r_V:A_*^G(X,i) \rightarrow A_*^G((X \times V^f),i)$$
and thus inverse systems $\{r'_V(K_i^{'G} \}_{ V \in {\cal V} }$
and $\{r_V(A_*^{'G}(X,i))\}_{ V \in {\cal V} }$.

\begin{thm} \label{completions}
There are isomorphisms of completions

$$\liminv r'_V(K_i^{'G}(X)) \simeq \widehat{K_i^{'G}(X)}$$
and if $X$ is quasi-projective we also have

$$\liminv r'_V(K_i^{'G}(X)) \simeq \widetilde{K_i^{'G}(X)}$$

$$\liminv r_V(A_*^{G}(X,i)) \simeq \widehat{A_*^{G}(X,i)}
\simeq \widetilde{A_*^{G}(X,i)}.$$
\end{thm}

Remark: A similar equality of completions was proved (for $K_0^{'G}$)
in \cite{CEPT} for actions of finite groups of
projective varieties defined over rings of integers of number fields.

\medskip

As a result of this identification of completions, we can prove
a particular case of a
conjecture of K\"ock (\cite{kock}) for arbitrary reductive groups acting
on regular schemes of finite type over a field. Set
$K(X,G) = \oplus K_i^G(X) = \oplus K_i^{'G}(X)$, and let
$\widetilde K(X,G)$ be the completion along the augmentation
ideal $\tilde{P}$ of $K_0^G(X)$.
\begin{cor} Let $X \stackrel{f} \rightarrow Y$ be a proper equivariant
morphism of quasi-projective, regular schemes. Then
$f_*:K(X,G) \rightarrow K(Y,G)$
is continuous with respect to the $\tilde{P}$-adic topologies.
\end{cor}
Proof. The pushforward
$f_*$ induces a map of inverse systems
$$\liminv r'_V(K_i^{'G}(X)) \rightarrow \liminv R'_V(K_i^{'G}(Y))$$
The corollary follows from the identification of completions in
Theorem \ref{completions}. \endproof

\medskip

Proof of Theorem \ref{completions}.
The two statements have essentially identical proofs, so we
will only prove the isomorphisms in $K$-theory. Furthermore,
the proof that
$$\liminv r'_V(K_i^{'G}(X)) \simeq \widetilde{K_i^{'G}(X)}$$
is virually identical to the proof that
$$\liminv r'_V(K_i^{'G}(X)) \simeq \widehat{K_i^{'G}(X)}$$
so we will only prove the latter. (The only difference in the proof
of the former is that we need to assume $X$ is quasi-projective
so we can compare the $\gamma$ filtration and the topological filtration
on $K_0^G(X)$.)
This statement is the analogue of
\cite[Theorem 2.1]{A-S}, except that we do not need the hypothesis
that $\widehat{K_i^{'G}(X)}$ is finite over $R(G)$.  As in \cite{A-S}
we first prove the result for a torus and from this deduce the general
case.  Our proof of the torus case is somewhat different from that of
\cite{A-S}, but the passage to the general case uses their arguments,
which we have repeated for completeness.

\medskip

{\em Step 1.} We first prove the result if $G=T$ is a torus.
We have filtrations of $K_i^{'G}(X)$ by the ideals $\mbox{ker }r'_V$
and by powers of the ideal $P = P_T$.  It suffices to show that the
filtrations have bounded difference.  This is a consequence of the
next two lemmas.

\begin{lemma}\label{koszul} Let $V$ be a representation of $T$ and $W
  \subset V$ a subrepresentation of codimension $l$.  Let $i: X \times
  W \rightarrow X \times V$ be the inclusion.  Then $i_*(K_i^{'T}({X
    \times W})) \in P^l K_i^{'T}(X \times V)$.
\end{lemma}
Proof of Lemma \ref{koszul}: We can find a chain of $T$-invariant
subspaces $W=W_l \subset W_{l-1} \subset \ldots \subset W_0 = V$ where
the codimension of $W_j$ in $V$ is $j$.  By induction on codimension,
it suffices to consider the case where the codimension of $W$ is $1$.
By the projection formula for equivariant $K$-theory, $i_*(K_i^{'T}(X
\times W))= i_*([{\cal O}_{X \times W}])K_i^{'T}(X)$.  However,
$$
i_*([{\cal O}_{X \times W}]) = [{\cal O}_{X \times V}] - [(V/W)
\otimes_{k} {\cal O}_{X \times V}],
$$
which is in $P K_i^{'T}(X \times V)$.  \endproof

\begin{lemma}\label{bound}
$$P^{s}K_0^{'T} \subset \mbox{ker } r'_V \subset P^l K_0^{'T}(X)$$
for any $s > d +  \mbox{dim }X - \mbox{dim }T$.
\end{lemma}

Proof of Lemma \ref{bound}.  If $V$ is a representation of a torus
then $V^u=V - V^f$ is a finite union of linear subspaces
(Appendix,  Proposition \ref{complement})
which by assumption have codimension at least $l$ in $V$. From the localization
long exact sequence
$$\ldots \rightarrow K_i^{'T}(X \times V^u) \stackrel{i_*} \rightarrow
K_i^{'T}(X \times V)
\stackrel{r'_V}\rightarrow K_i^{'T}(X \times V^f) \rightarrow \ldots$$
we know
that $\mbox{ker }r'_V\; = i_*(K_i^{'T}(X \times V^u))$. The image of
$K_i^{'T}(X \times V^u)$ is generated by the images of $K_i^{'T}(X \times W)$
for each linear space $W \subset V^u$.  By Lemma \ref{koszul}, these
images are contained in $P^l K_i^{'T}(X \times V) = K_i^{'T}(X)$.
Hence $\mbox{ker }r'_V \subset P^l K_i^{'T}(X)$.

For the other inclusion, note that $K_i^{'T}(X \times V^f) =
K_i'(X_T)$, and $P^s K_0^{'T}(X \times V^f) \subset F^s K_i'(X_T)$,
where $F^{\cdot}$ denotes the $\gamma$-filtration on $R(G)$.
Since a point is projective,
$F^sK_i'(X_T) \subset F^s_{top}K_i^{'T}(X_T)$. Thus, if
$s > \mbox{dim }X_T$, then  $F^s K_i'(X_T) = 0$.  Hence
$P^{s}K_i^{'T} \subset \mbox{ker } r'_V$, as desired. \endproof

This lemma implies the desired equality of completions for the case of
a torus.

\medskip

{\em Step 2.}  We prove the result for $G=GL_n$.  Let $j: T
\hookrightarrow G$ be the inclusion and let $j^*: K_i^{'G}(X)
\rightarrow K_i^{'T}(X)$ be the induced restriction map.

\begin{lemma} \label{summand}
There is a functorial map $j_*: K_i^{'T}(X) \rightarrow K_i^{'G}(X)$
such that $j_*j^*$ is the identity.  Hence $K_i^{'G}(X)$ is a direct
summand in $K_i^{'T}(X)$.
\end{lemma}

Proof of Lemma \ref{summand}.
This is proved in \cite[Prop. 4.9]{Atiyah} for topological K-theory.
The same proof works in the algebraic setting: the main ingredient is the
projective bundle theorem, which was proved in this setting by Thomason
\cite[Theorem 3.1]{Thomason}. \endproof

{}From the proof of Step 1, in computing $\liminv K_i^{'T}(X \times
V^f)$ we need not consider all representations of $T$: it suffices to
consider the subsystem of representations of $T$ which are
restrictions of representations of $G$,
then $\mbox{ker }r^{'G}_V = \mbox{ker }r^{'T}_V \cap K_i^{'G}(X)$.

The submodule $K_i^{'G}(X)$ of $K_i^{'T}(X)$ inherits two topologies
from $K_i^{'T}(X)$: the topology induced by the ideals $\mbox{ker
  }r^{'T}_V \cap K_i^{'G}(X) = \mbox{ker }r^{'G}_V$, and the topology
induced by powers of the ideal $P_T$.  Because $K_i^{'G}(X)$ is a
direct summand in $K_i^{'T}(X)$, by Lemma \ref{bound} these topologies
coincide.  On the other hand, as noted in \cite{A-S}, the ideals $P_T$
and $P_G R(T)$ have bounded difference,
so they induce the
same topology on $K_i^{'T}(X)$.  The restriction of this topology to
$K_i^{'G}(X)$ is the topology induced by powers of the ideal $P_G$.
Putting these facts together, we conclude that $\liminv K_i^{'G}(X
\times V^f) \simeq \widehat{K_i^{'G}(X)}$ for $G=GL_n$, as desired.

\medskip

{\em Step 3.}  We now deduce the result for general $G$.  Embed $G$
into $H = GL_n$.  Then $K_i^{'G}(X) = K_i^{'H}(X \times^G H)$\footnote{Note
that $X \times^G H$ is a
scheme, because of our hypothesis on $X$ or $G$.}
(\cite[Proposition 6.2]{Thomason}).

As above, we may restrict our attention
to representations of $G$ which are restrictions of representations of
$H$, then
$$\liminv r'_V(K_i^{'G}(X)) \simeq \liminv r'_V(K_i^{'H}(( X \times^G H)).$$
As noted in \cite{A-S}, the $P_H$-adic and $P_G$-adic topologies
coincide on any $R(G)$-module, and hence, by the result for $H=GL_n$,
we have
$$\liminv r'_V(K_i^{'G}(X)) \simeq \widehat{K_i^{'G}(X)},$$
as desired.  \endproof

If $X$ is any $G$-scheme, let $K_i'^G(X)_{\Q} =
K_i'^G(X) \otimes \Q$, and $A_*^G(X,i)_{\Q} = A_*^G(X,i) \otimes \Q$.

\begin{cor} There are isomorphisms of completions
$$\liminv \{ r'_V(K_i^{'G}(X))_{\Q} \} \simeq \widehat{K_i^{'G}(X)_{\Q}}
  \simeq \widetilde{K_i^{'G}(X)_{\Q}}$$
and if $X$ is quasi-projective
$$\liminv \{ r_V(A_*^{G}(X,i))\} \simeq
\widehat{A_*^{G}(X,i)_{\Q} }.$$
\end{cor}
Proof: The proof is the same as above, except we do not need $X$
to be quasi-projective to show that
$\tilde{P}^sK_i^{'G}(X)_{\Q} \subset \mbox{ker }r'_V$ for $s >> 0$. Instead,
we can apply the non-equivariant Riemann-Roch isomorphism of
\cite[Chapter 18]{Fulton} \endproof

\medskip

\noindent{\bf Example.}
Since inverse limits do not commute with tensoring with $\Q$,
$\liminv \{ r'_V(K_i^{'G}(X))_{\Q}
\}$ need not be equal to $ (\liminv r'_V(K_i^{'G}(X))) \otimes
\Q$. For example if  $G = \Z / 2\Z$ and $X = pt$
then
$$K_0^{'G}(X)  = K_0^G(X) = R(G) = \Z[u]/(u^2-1).$$
In this case $\liminv r'_V(K_0^{'G}(X)) \simeq \Z_{(2)}$ (the
2-adic integers). Thus, $(\liminv r'_V(K_0^{'G}(X))) \otimes \Q = \Q_{2}$.
On the other hand, $\liminv(r'_V(K_0^{'G}(X)_{\Q})= \Q$.

However, there is a map
$$
\liminv r'_V(K_i^{'G}(X)) \rightarrow \liminv \{ r'_V(K_i^{'G}(X))_{\Q} \}
$$
which induces a map
$$
( \liminv r'_V(K_i^{'G}(X))) \otimes \Q \rightarrow \liminv \{
r'_V(K_i^{'G}(X))_{\Q} \}.
$$
Likewise there is a map
$$
( \widehat{K_i^{'G}(X)} ) \otimes \Q \rightarrow\widehat{K_i^{'G}(X)_{\Q}}
$$
which need not be an isomorphism.
\subsection{The equivariant Riemann-Roch isomorphism} \label{eqrr}
Before stating the Riemann-Roch theorem we need to define the equivariant
Chern character.  For this we need a suitable completion of
$A^*_G(X)$.  Elements of this completion should operate on
$\widehat{A_*^G(X)}$.  There are three candidates for this completion:

(1) $\liminv A^*_G(X \times V^f)$;

(2) $\widetilde{A^*_G(X)}$, the completion with respect to $\tilde{Q}$;

(3) $\widehat{A^*_G(X)}$, the completion with respect to $Q$. \\
If $X$ is smooth, then all three completions are equal. In general,
we do
not know whether they are equal because of the lack
of a suitable exact sequence for operational Chow groups.  However,
all of these operate on $\widehat{A_*^G(X)}$ by virtue of the
isomorphisms proved in the last subsection.  There are maps
$$
\widehat{A^*_G(X)} \rightarrow \widetilde{A^*_G(X)} \rightarrow
\liminv A^*_G(X \times V^f).
$$
The first map is because $Q \subset \tilde{Q}$.  The second is because
the map $A^*_G(X) \rightarrow A^*_G(X \times V^f)$ induces a map
$\widetilde{A^*_G(X)} \rightarrow A^*_G(X \times V^f)$ (because high
powers of $\tilde{Q}$ map to zero in $A^*_G(X \times V^f)$).  These
maps are compatible with restrictions and hence induce a map to the
inverse limit.  We will define a Chern character map with image in
$\widetilde{A^*_G(X)}$.

\begin{defn}
Define the equivariant Chern character
$$
ch_G:K_0^G(X) \rightarrow \widetilde{A^*_G(X)}_{\Q}$$
by the formula
$$ch_G(E) = r + c_1^G(E) + \frac{1}{2}(c_1^G(E)^2 - 2c_2^G(E))
+ \ldots.$$

\end{defn}

Let $V$ be a representation of $G$ such that $V -V^f$ has codimension
more than $i$. If $E \rightarrow X$ is an equivariant vector bundle,
then $c_i^G(E)$ restricts to $c_i(E \times^G  V^f)$
under the restriction map $A_G^i(X) \rightarrow A_G^i(X \times^G V^f)$.
Thus, the equivariant Chern character $ch_G:K_0^G(X) \rightarrow \widetilde{
A_G^*(X)}$ restricts to the ordinary Chern character
$ch_{X \times^G V^f}:K_0(X \times^G V^f) \rightarrow A^*(X \times^G V^f)$.

\begin{prop}

There is a factorization
$$ch_G:K_0^G(X) \rightarrow \widehat{K_0^G(X)} \rightarrow
\widetilde{K_0^G(X)} \widetilde{A^*_G(X)}_{\Q}.$$
\end{prop}
Proof. The proof follows from the fact that $ch(P^n)$ and
$ch(\tilde{P}^n)$ are contained in $\tilde{Q}^n$ for any $n>0$.
\endproof

We will also denote the map $\widehat{K_0^G(X)} \rightarrow
\widetilde{A^*_G(X)}_{\Q}$ by $ch_G$.

\begin{thm} \label{rockandroll} (Equivariant Riemann-Roch)\\

There are maps
$$\t^G_X: \widehat{K_{0}^{'G}(X)}
\rightarrow \widehat{A_*^G(X)_{\Q}}$$
with the following properties (cf. \cite[Chapter 18]{Fulton}):

(1) $\t^G_X$ is covariant for equivariant proper morphisms.

(2) If $\epsilon \in K_0^G(X)$ and $\alpha \in \widehat{K_0^{'G}(X)}$
then $\t^G_X(\epsilon \alpha) = ch^G_X(\epsilon) \cap \t_X^G(\alpha)$.
(Recall that $\widetilde{A_G^*(X)}$ operates on $\widehat{A_*^G(X)}$
because of the isomorphisms of completions.)

(3) If $f:X \rightarrow Y$ is a $G$-equivariant l.c.i. morphism, then
$$ch^G f_*(\epsilon) = f_*(((Td^G(T_f) ch^G(\epsilon))$$
and
$$\t_X^G(f^*\alpha) = Td^G(T_f)f^*\t_X^G(\alpha).$$

(4) If $V \subset X$ is a $G$-invariant subvariety of dimension $k$, then
$$\t_X^G({\cal O}_V - [V]_G)  \in F_{k-1}(\widehat{A_*^G(X)_{\Q}})$$
where $F_{k-1}(\widehat{A_*^G(X)_{\Q}})$ denotes the subgroup of cycles
of ``dimension'' strictly less than $k$.

(5) $\t_X^G$ factors through the map $\widehat{K_0^{'G}(X)}
\rightarrow\widehat{K_0^{'G}(X)_{\Q}}$ and induces an isomorphism
between $\widehat{K_0^{'G}(X)_\Q}$ and $\widehat{A_*^G(X)_\Q}$.

(6) If $X$ is quasi-projective, and $i > 0$, then there is an
isomorphism
$$\t_X^G: \liminv{K_i^{'G}(X \times V^f)_\Q} \stackrel{\simeq} \rightarrow
\liminv{A_*^G((X \times V^f),i)_\Q}.$$
\end{thm}

Proof. When $i = 0$, the restriction maps $r_V$ and $r'_V$ are surjective.
Thus
by Theorem \ref{completions},
$$\widehat{K_0^{'G}(X)} = \liminv K_0{'G}(X \times V^f)$$
and
$$\widehat{A_*^G(X)_{\Q}} = \liminv (A_*^G(X \times V^f) \otimes \Q).$$
By non-equivariant Riemann-Roch (\cite[Chapter 18]{Fulton}),
for each  representation $V$ of $G$ there is a
map
$$\tau_{X \times V_f}: K_0^{'G}(X \times V^f) \rightarrow A_*^G(X
\times V^f) \otimes \Q$$
satisfying the analogues of (1) - (5). To prove the theorem it
suffices to show that maps $\{\tau_{X \times V^f}\}$ are compatible
with the inverse system maps.  This is quite straightforward.  Let $V$
and $W$ be representations. Let $\pi:(X \times V^f \oplus W)/G
\rightarrow (X \times V^f)$. Then $\pi$ is smooth and $Td_{T_{\pi}} =
1$. Thus $\pi^* \cdot (\tau_{X \times V_f}) = \tau_{X \times V^f
  \oplus W} \cdot \pi^*$. Likewise, if $i:(X \times^G V^f \oplus W)
\rightarrow (X \times^G (V \oplus W)^f)$ is the inclusion map, then
$i^* \cdot \tau_{X \times (V \oplus W)^f} =\tau_{X \times V^f \oplus
  W} \cdot i^*$.

To prove (6) we argue as above, using Bloch's Riemann-Roch isomorphism
$$K_i(X \times^G V^f) \cong A_*((X \times^G V^f),i).$$
\endproof

\paragraph{Vistoli's conjecture}
By composing the map above with the natural map $K_0^{'G}(X)
\rightarrow \widehat{K_0^{'G}(X)}$ we get a map
$\t^G_X: K_{0}^{'G}(X)
\rightarrow \widehat{A_*^G(X)_{\Q}}.$
When $G$ acts on $X$ with finite reduced stabilizers then Vistoli
\cite{Vi3} stated a theorem which asserted the existence of
a map
$$\t_X: K_0^{'G}(X) \otimes \Q \rightarrow A_*([X/G] \otimes \Q)$$
satisfying properties (1)-(4) above (here $[X/G]$ is the
Deligne-Mumford quotient stack).  By Theorem \ref{moduli}, $A_*([X/G])
\otimes \Q= A_*^G(X) \otimes \Q$. Thus $I^d(A_*^G(X) \otimes \Q)=0$
for $d >>0$, so $\widehat{A_*^G(X)_{\Q}} = A_*^G(X) \otimes \Q$. Thus
Vistoli's map is a special case of our map $\tau_X^G$, since it is
uniquely determined by properties (1)-(4).  Vistoli noted that this
map need not be an isomorphism and made the following conjecture
about its kernel.

\begin{conj} (\cite[Conjecture 2.4]{Vi3})
Suppose that $G$ acts on $X$ with finite reduced stabilizers.
If $\alpha \in ker(\tau^G_X: K_0^{'G}(X) \rightarrow A_*([X/G])_{\Q})$
then there exists an element $\epsilon \in K_0^G(X)$
with every non-zero rank (meaning $\epsilon$ is represented by
a complex of locally free sheaves whose homology is non-zero
at the generic point of every subvariety) such that $\epsilon \alpha = 0$.
\end{conj}

The results of this section identify
the kernel: it is exactly the kernel of the completion map
$K_0^{'G}(X)\otimes \Q \rightarrow \widehat{K_0^{'G}(X)_\Q}$.

\begin{prop}
Suppose $K_0^G(X)$ is Noetherian and $K_0^{'G}(X)$
is finitely generated over $K_0^{G}(X)$.  Then
$\alpha \in ker\; \tau_X^G$ if and only if $(1+\delta) \alpha = 0$
for some $\delta \in K_0^G(X)$ of (virtual) rank 0.
\end{prop}
Proof.
The proof follows immediately from Krull's theorem. \endproof

\section{Localization}
In this section we discuss properties of equivariant Chow groups
that are similar to properties of equivariant cohomology.
In the first part, we give the relationship between
$A_*^G(X)$ and $A_*^T(X)$ when $G$ is a connected reductive
group with maximal torus $T$. The remainder of the section
is devoted to actions of (split) tori. In particular, we prove
two localization theorems (Theorems \ref{lcztn}, \ref{xxx}).
Following ideas of \cite{A-B} they
yield a characteristic free proof of the Bott residue
formula for split torus actions on complete varieties over a field
of arbitrary characteristic (Theorem \ref{bott}).

\subsection{Connected reductive groups}
Denote by $A^*_G$ or $R_G$ the equivariant Chow ring of a point (the
equivariant Chow groups of a point have a ring structure since a point
is smooth). If $G$ is a connected reductive group then by
\cite{E-G}, $R_G \otimes \Q=
Sym(\hat{T})^W \otimes \Q$, where $\hat{T}$ is the group of characters
of the maximal torus and $W$ is the Weyl group. When $G$ is special
in the sense of \cite{Sem-Chev} then $R_G = Sym(\hat{T})^W$ exactly
(\cite{E-G}). Under this identification we will write
$R_G^d = A^d_G = A^G_{-d}$.  Via pullback
from a point, $A_*^G(X)$ has the structure of an $R_G$-module.

If $G = T$ is a split torus, then $W$ is trivial, and the identification
$R_T = Sym(\hat{T})$ is given explicitly as follows.  If $\lambda \in \hat{T}$,
let $k_{\lambda}$ denote the corresponding 1-dimensional
representation of $T$, and let $L_{\lambda}$ denote the line bundle
$U \times^T k_{\lambda} \rightarrow U /T$.   The map $\hat{T}
\rightarrow A^1_T$ given by $\lambda \mapsto c_1(L_{\lambda})$ extends
to a ring isomorphism $Sym(\hat{T}) \rightarrow R_T$.  If $f:T
\rightarrow S$ is a homomorphism of tori, then there is a pullback map
$f^*:\hat{S} \rightarrow \hat{T}$.  This extends to a ring homomorphism
$f^*: Sym(\hat{S}) \rightarrow Sym(\hat{T})$, or in other words, a map
$f^*: R_S \rightarrow R_T$.

\begin{prop} Let $G$ be a connected reductive group with split maximal
torus $T$ and Weyl group $W$. Then
$A_*^G(X) \otimes \Q = A_*^T(X)^W \otimes \Q$. If $G$ is special
the isomorphism holds with integer coefficients.
\end{prop}

Proof: If $G$ acts freely on $U$, then so does
$T$. Thus for a sufficiently large representation $V$,
$A_{i}^T(X) = A_{i+l -t}((X \times U)/T)$ and
$A_i^G(X) = A_{i+l-g}((X \times U)/G)$ (here
$l$ is the dimension of $V$, $t$ the dimension of $T$
and $g$ the dimension of $G$). On the other hand,
$(X \times U)/T$ is $G/T$ bundle over $(X \times U)/G$.
Thus
$A_{k}((X \times U/T)) \otimes \Q =A_{k+g-t}((X \times U)/G)^W \otimes \Q$
and if $G$ is special, then the equality holds integrally (\cite{E-G})
and the proposition follows.
\endproof

\medskip

Thus, for connected reductive groups,
to compute equivariant Chow groups (at least
with rational coefficients), it suffices to understand
equivariant Chow groups for tori. We begin with the following
proposition.

\begin{prop}
If $T$ acts trivially on $X$, then $A_*^T(X) = A_*(X) \otimes
R_T$.
\end{prop}

Proof. If the action is trivial then $(U \times X)/T= U/T \times X$.
The spaces $U/T$ can be taken to be products
of projective spaces, so $A_*(U/T \times X) = A_*(X) \otimes A_*(U/T)$.
\endproof

\medskip

{\bf Remark.} If the action is trivial,
the pullbacks $A^*X_T \rightarrow A^*X$ and $A^*(U/T) \rightarrow
A^*X_T$ induce an inclusion of $A^*X \otimes R_T
\subset A^*_T(X)$ as a subring. If $X$ is smooth, then
the inclusion is an isomorphism by Proposition \ref{opsmooth}.

\subsection{Fixed loci and the localization theorem}
For the remainder of this section,
all Chow groups have rational coefficients, and
for simplicity of exposition, we assume that tori are split.

If $X$ is a scheme with a $T$-action,
we may put a closed subscheme structure
on the locus $X^T$ of points fixed by
$T$.

Now $R_T= Sym(\hat{T})$ is a polynomial ring.
Set ${\cal Q}= (R_T^+)^{-1} \cdot R_T$, where $R_T^+$
is the multiplicative system of homogeneous elements of positive degree.
\begin{thm} \label{lcztn}(localization)
The map $i^T_*:A_*(X^T)  \otimes {\cal Q} \rightarrow A_*^T(X) \otimes
{\cal Q}$ is
surjective, and if $X$ is quasi-projective it is an isomorphism.
\end{thm}

{\bf Remark.}  The quasi-projectivity assumption is needed to apply
the long exact sequence for higher Chow groups.  The strategy of the
proof is similar to \cite[Theorem 5.3]{Th2}.

\medskip

Proof. Applying the localization exact sequence for higher equivariant
Chow groups
$$\ldots \rightarrow A_*^T(X^T) \rightarrow A_*^T(X) \rightarrow
A_*^T(X-X^T) \rightarrow 0$$
the theorem follows from the following proposition.

\begin{prop} \label{fix}
If $T$ acts on $X$
without fixed points, then there exists $r \in R_T^+$ such that
$r \cdot A_*^T(X,m)= 0$. (Recall that $A_*^T(X,m)$ refers to
$T$-equivariant higher Chow groups.)
\end{prop}

Suppose $f: T \rightarrow S$ is a homomorphism of tori.  As discussed
above, there is a pullback map $f^*: A^*_S \rightarrow A^*_T$.

\begin{lemma} \label{t-map}(cf. \cite{A-B})
Suppose there is a $T$-map
$X \stackrel{\phi} \rightarrow S$.
Then $t\cdot A_*^T(X)= 0$ for any $t=f^*s$ with $s \in R_S^+$.
\end{lemma}

Proof of Lemma \ref{t-map}.
Since $A^*_S$ is generated in degree 1, we may
assume that $s$ has degree 1. After clearing denominators
we may assume that $s = c_1(L_s)$ for some line bundle
on a space $U/S$. The action of  $t=f^*s$ on $A_*(X_T)$ is just
given by $c_1(\pi_T^*f^*L_s)$ where $\pi_T$ is the
map $U \times^T X \rightarrow U/T$.
To prove the lemma we will show that this bundle is trivial.

First note that  $L_s = U \times^S k$ for some action of $S$ on
the one-dimensional vector space $k$.
The pullback bundle on $X_T$ is the line bundle
$$U \times^T(X \times k) \rightarrow X_T$$
where $T$ acts on $k$ by the composition of $f:T \rightarrow S$
with the original $S$-action.
Now define a map
$$\Phi: X_T \times k \rightarrow U \times^T(X \times k)$$
by the formula
$$\Phi(e,x,v) = (e,x,\phi(x)\cdot v)$$
(where $\phi(x) \cdot v$ indicates the original $S$ action).
This map is well defined since
\begin{eqnarray*}
\Phi(et,t^{-1}x,v) & = &(et, t^{-1}x, \phi(t^{-1}x) \cdot v)\\
& = & (et,t^{-1}x,t^{-1} \cdot(\phi(x) \cdot v))
\end{eqnarray*}
as required. This map is easily seen to be an isomorphism
with inverse $(e,x,v) \mapsto (e,x,\phi(x)^{-1} \cdot v)$.
\endproof

\medskip

Proof of Proposition \ref{fix}. Since $A_*^G(X) = A_*^G(X_{red})$
we may assume $X$ is reduced. Working with each component
separately, we may assume $X$ is a variety. Let $X^0 \subset X$
be the ($G$-invariant) locus of smooth points.
By Sumhiro's theorem \cite{Sumihiro}, the action of
a torus on a normal variety is locally
linearizable (i.e. every point has an affine invariant
neighborhood). Using this theorem it is easy to see
that the set
$X(T_1) \subset X^0$ of points with stabilizer
$T_1$ can be given the structure of a locally closed
subscheme of $X$.
Furthermore, only finitely many subgroups can occur as stabilizers
(Appendix, Lemma \ref{porb}), so there is
some $T_1$ such that $U= X(T_1)$ is open in $X^0$, and thus in $X$.

The torus $T'=T/T_1$ acts without stabilizers,
but the action of $T'$ on $U$ is not a priori proper. However,
by \cite[Proposition 4.10]{Th1}, we can replace $U$
by a sufficiently small open set so that $T'$ acts freely
on $U$ and a principal bundle quotient $U \rightarrow U/T$
exists. Shrinking $U$ further, we can assume that this
bundle is trivial, so there is a $T$ map $U \rightarrow T'$.
Hence, by the lemma, $t \cdot A^T_*(U) = 0$ for any $t \in A_T^*$
which is pulled back
from $A^*_{T'}$.

Let $Z = X -U$.  By induction on dimension, we may assume $p \cdot
A^T_*Z = 0$ for some homogeneous polynomial $p \in R_T$.  From the
long exact sequence of higher Chow groups,
$$\ldots A_*^T(Z,m) \rightarrow A_*^T(X,m) \rightarrow A_*^T(U,m) \rightarrow
\ldots$$
it follows that $tp$ annihilates $A_*^T(X)$ where
$t$ is the pullback of a homogeneous element of degree $1$
in $R_S$.

\medskip

Remark: Using only the short exact localization sequence for ordinary
equivariant Chow groups (which does not require an assumption
of quasi-projectivity) shows that $i_*$ is surjective.
\endproof

\subsection{Explicit localization and the integration formula}
The localization theorem in equivariant cohomology has a more explicit
version for manifolds.  This yields an integration formula from which
the Bott residue formula is easily deduced (\cite{A-B}, \cite{B-V}).
In this section we prove the analogous results for equivariant Chow
groups of smooth varieties.  Because equivariant Chow theory has
formal properties similar to equivariant cohomology, the arguments are
almost the same as in \cite{A-B}. As before we assume that all tori
are split.

Let $F$ be a scheme with a trivial $T$-action.
If $E \rightarrow F$ is a $T$-equivariant vector bundle on
$F$, then $E$ splits canonically into a direct sum of vector subbundles
$\oplus_{\lambda \in \hat{T}} E_{\lambda}$, where $E_{\lambda}$
consists of the subbundle of vectors in $E$ on which $T$ acts by the
character $\lambda$.  The equivariant Chern classes of an
eigenbundle $E_{\lambda}$ are given by the following lemma.

\begin{lemma} \label{l.trivchern}
Let $F$ be a scheme with a trivial $T$-action, and let
$E_{\lambda} \rightarrow F$ be a $T$-equivariant vector bundle of rank
$r$ such
that the action of $T$ on each vector in $E_{\lambda}$ is given by the
character $\lambda$.  Then for any $i$,
$$
c^T_i(E_{\lambda}) = \sum_{j \leq i}
\left( \begin{array}{c} r-j \\
i-j \end{array}
\right)
c_j(E_{\lambda}) \lambda^{i-j}.
$$
In particular the component of $c_r^T(E_{\lambda})$ in $R^r_T$ is
given by $\lambda^r$.  \endproof
\end{lemma}

As noted above, $A^*_T(F) \supset A^*F \otimes R_T$.  The lemma
implies that $c^T_i(E)$ lies in the subring $A^*F \otimes R_T$.
Because $A^N F = 0$ for $N > \mbox{dim }F$, elements of $A^i F$, for
$i>0$, are nilpotent elements in the ring $A^*_T(F)$.  Hence an
element $\alpha \in A^d F \otimes R_T$ is invertible in $A^*_T(F)$ if
its component in $A^0 F \otimes R^d_T \cong R^d_T$ is nonzero.

For the remainder of this section $X$ will denote a smooth variety
with a $T$ action.  If $X$ is smooth then by \cite{Iv} the fixed locus
$X^T$ is also smooth.  For each component $F$ of the fixed locus $X^T$
the normal bundle $N_FX$ is a $T$-equivariant vector bundle over $F$.
Note that the action of $T$ on $N_FX$ is non-trivial.

\begin{prop}
If $F$ is a component of $X^T$
with codimension $d$ then $c_d^T(N_FX)$ is invertible
in $A^*_T(F) \otimes {\cal Q}$.
\end{prop}
Proof: By (\cite[Proof of Proposition 1.3]
{Iv}), for each closed point $f \in F$, the tangent space
$T_fF$ is equal to $(T_fX)^T$, so $T$ acts with non-zero weights on the
normal space $N_f = T_fX/T_fF$. Hence the characters $\lambda_i$
occurring in the
eigenbundle decomposition of $N_FX$ are all non-zero.  By the
preceding lemma, the component of $c_d^T(N_FX)$ in $R^d_T$ is nonzero.
Hence $c_d^T(N_FX)$ is invertible
in $A^*_T(F) \otimes {\cal Q}$, as desired.  \endproof

\medskip

Using this result we can get, for $X$ smooth,
the following more explicit version
of the localization theorem.

\begin{thm} \label{xxx}(Explicit localization)
Let $X$ be a smooth (not necessarily quasi-projective) variety with a torus
action.
Let $\alpha \in A_*^T(X) \otimes {\cal Q}$.
Then $$\alpha = \sum_F
i_{F*}\frac{i^*_F\alpha}{c_{d_F}^T(N_FX)},$$ where the sum is over the
components $F$ of $X^T$ and $d_F$ is the codimension of $F$ in $X$.
\end{thm}
Proof: By the surjectivity part of the localization theorem,
we can write $\alpha = \sum_F
i_{F*}(\beta_F)$.  Therefore, $i^*_F\alpha = i^*_Fi_{F*}(\beta_F)$
(the other components of $X^T$ do not contribute); by the
self-intersection formula, this is equal to $ c_{d_F}^T(N_FX) \cdot
\beta_F$.  Hence $\beta_F = \frac{i^*_F\alpha}{c_{d_F}^T(N_FX)}$ as
desired. \endproof

\medskip

If $X$ is complete, then the projection $\pi_X: X \rightarrow pt$
induces push-forward maps $\pi^T_{X*}: A^T_* X \rightarrow R_T$ and
$\pi^T_{X*}: A^T_* X \otimes {\cal Q} \rightarrow {\cal Q}$.  There
are similar maps with $X$ replaced with any component $F$ of $X^T$.
Applying $\pi^T_{X*}$ to both sides of the explicit localization
theorem, and noting that $\pi^T_{X*} i_{F*} = \pi^T_{F*} $, we deduce
the ``integration formula'' (cf. \cite[Equation (3.8)]{A-B}).

\begin{cor}
(Integration formula) Let $X$ be smooth and complete, and
let $\alpha \in A_*^T(X) \otimes {\cal Q}$.  Then
$$\pi_{X*}(\alpha) = \sum_{F \subset X^T}
\pi_{F*}\{\frac{i^*_F\alpha}{c_{d_F}^T
(N_FX)}\}$$
as elements of ${\cal Q}$.  \endproof
\end{cor}

\medskip

{\bf Remark.} If $\alpha$ is in the image of the natural map $A_*^T(X)
\rightarrow A_*^T(X) \otimes {\cal Q}$ (which need not be injective),
then the equation above holds in the subring $R_T$ of ${\cal Q}$.  The
reason is that the left side actually
lies in the subring $R_T$; hence so does the right side.  In the
results that follow, we will have expressions of the form $z = \sum
z_j$, where the $z_j$ are degree zero elements of ${\cal Q}$ whose sum
$z$ lies in the subring $R_T$.  The pullback map from equivariant to
ordinary Chow groups gives a map $i^*: R_T = A^T_* (pt) \rightarrow \Q
= A_* (pt)$, which identifies the degree 0 part of $R_T$ with $\Q$.
Since $\sum z_j$ is a degree 0 element of $R_T$, it is identified via
$i^*$ with a rational number.  Note that $i^*$ cannot be applied to
each $z_j$ separately, but only to their sum.  In the integration and
residue formulas below we will identify the degree 0 part of $R_T$
with $\Q$ and suppress the map $i^*$.  \medskip

The preceding corollary yields an integration formula for
an element $a$ of the ordinary Chow group $A_0 X$, provided that $a$ is
the pullback of an element $\alpha \in A^T_0 X$.

\begin{prop}
Let $a \in A_0 X$, and suppose that $a = i^* \alpha$ for $\alpha \in
A^T_0 X$.  Then
$$
\mbox{deg }(a) = \sum_F  \pi^T_{F*}\{\frac{i^*_F\alpha}{c_{d_F}^T
(N_FX)} \}
$$
\end{prop}
Proof: Consider the commutative diagram
$$\begin{array}{ccc}
X & \stackrel{i} \hookrightarrow & X_T\\
\downarrow\scriptsize{\pi_X} & & \downarrow\scriptsize{\pi^T_X}\\
\mbox{pt} & \stackrel{i} \rightarrow & U/T .
\end{array}$$
We have $\pi_{X*}(a) = \pi_{X*} i^*(\alpha) = i^* \pi^T_{X*}(\alpha)$.
Applying the integration formula gives the result. \endproof

\subsection{The residue formula}
Let $E \rightarrow X$ be a $T$-equivariant vector bundle of rank $r$
on a complete, smooth $n$-dimensional variety. Let $p(x_1, x_2, \ldots
, x_r)$ be a polynomial of weighted degree $n$, where the degree of
$x_i$ is $i$.  The integration formula above will allow us to compute
$\mbox{deg }(p(c_1(E), \ldots, c_r(E)) \cap [X])$
in terms of the
restriction of $E$ to $X^T$.

As a notational shorthand, write $p(E)$
for $p(c_1(E), \ldots, c_r(E))$ and $p^T(E)$ for $p(c^T_1(E), \ldots,
c^T_r(E))$.
Write $p^T(E|_{F})$ for $p(c^T_1(E|_{F}), \ldots,
c^T_r(E|_{F})) = i^*_F p(c^T_1(E), \ldots, c^T_r(E))$.  Notice that
$p(E) \cap [X] = i^* (p^T(E) \cap [X_T])$.  We can therefore apply the
preceding proposition to get the Bott residue formula.

\begin{thm} \label{bott}
(Bott residue formula) Let $E \rightarrow X$ be a $T$-equivariant
vector bundle of rank $r$ on a complete, smooth $n$-dimensional
variety, and let $p(x_1, x_2, \ldots , x_r)$ be a polynomial of weighted
degree $n$.  Then
$$
\mbox{deg }(p(E) \cap [X]) =  \sum_{F \subset X^T}
\pi^T_{F*}\{\frac{p^T(E|_{F}) \cap
[F]_T}{c_{d_F}^T (N_FX)} \}  .
$$
\endproof
\end{thm}

By Lemma \ref{l.trivchern} the
equivariant Chern classes $c^T_i(E|_{F})$ and $c_{d_F}^T (N_FX)$ can
be computed in terms of the characters of the
torus occurring in the eigenbundle
decompositions of $E|_{F}$ and $N_FX$ and the Chern classes of the
eigenbundles.
The above formula can then be readily converted (cf. \cite{A-B}) to more
familiar forms of the Bott residue formula not involving equivariant
cohomology. We omit the details. If the torus $T$ is
1-dimensional, then degree zero elements of ${\cal Q}$ are rational
numbers, and the right hand side of the formula is just a sum of
rational numbers. This is the form of the Bott residue formula which
is most familiar in practice.

\section{Actions of group schemes over an arbitrary base} \label{mixed}
Let $S$ be a regular scheme, and let $X/S$ be an $S$-scheme
with an action of a connected reductive group scheme $G/S$. With appropriate
assumptions on $G/S$ (see below), it is possible to use
Seshadri's results on geometric reductivity over an arbitrary
base to extend much of our theory.

If $S = Spec(\Z)$ and $G/S$ is reductive, then the
theory goes through more or less intact. In particular, if
$X/\Z$ is a smooth scheme acted on by a reductive group scheme
$G/\Z$, then there is an equivariant Chow ring $A^*_G(X)$. Such a ring
should be useful for studying intersection theory on moduli in
mixed characteristic.

\subsection{Definitions}

\begin{defn} Let $G/S$ be a smooth group scheme. Let $E/S$ be a vector
bundle (i.e., $Spec(\mbox{Sym}({\cal E^{.}})$
where ${\cal E}/S$ is a locally free
$G$-module). The bundle $E/S$ is said to be a representation of
$G/S$ if there is an action $G \times E \rightarrow E$ which
is linear on each fiber.
\end{defn}

We assume the following condition on $G/S$:

(*) There exist representations $E/S$ with a non-empty open
set $U/S$ such that $G/S$ acts freely on $U$.

\begin{prop}
If $G/S$ is a smooth group scheme, then condition (*) above is satisfied
if either

(1) $G/S$ is the pullback of
a group scheme $G_R/Spec(R)$ where $R$ is a Dedekind domain.

(2) The geometric fibers of $G/S$ are all semisimple with trivial
center.
\end{prop}

Proof. By \cite[Lemma 1, Proposition 3]{Seshadri}, $R(G_R)$ is
a projective $R$-module which contains a finitely generated $G$-invariant
$R$-module. Since $R$ is a Dedekind domain, this module is projective.
Pulling back to $S$ gives the desired representation.

By \cite[Expose II]{SGA3}, the Lie algebra $Lie(G/S)$ is a vector
bundle over $S$. Since $G$ has trivial center, the adjoint action of
$G$ on the vector bundle $Lie(G/S) \times_S Lie(G/S)$ is generically
free.
\endproof

\medskip

Henceforth we will assume that $G/S$ is reductive.
If condition (*) holds, we can find representations $E/S$ and
open set $U/S$ so that $E-U$ has arbitrarily high codimension.
By Seshadri's theorem (\cite{Seshadri}) there is a principal
bundle quotient $U \rightarrow U/G$.

The arguments of Proposition \ref{inap} yield
\begin{prop} \label{qbase}
Assume one of the following:

(1) $X/S$ equivariantly embeds in a projective bundle over $S$.

(2) $X/S$ is normal.

\noindent
Then a principal bundle quotient $X \times U \rightarrow X \times^G U$
exists.
\endproof \end{prop}

As a consequence of Proposition \ref{qbase} we can define
equivariant Chow groups.
\begin{defn} Assume that condition (*) on $G/S$ holds, as well
as one of the hypotheses (1) or (2) of Proposition \ref{qbase}.
Define the $i$-th equivariant Chow group as $A_{i+l-g}(X \times^G U)$
where $l = \mbox{dim }(U/S)$ and $g=\mbox{dim }(G/S)$. As for
algebraic schemes, the definition is independent of the representation.
\end{defn}

\subsection{Results over an arbitrary base}

Since most of the results of intersection theory hold for schemes
over a regular base (\cite[Chapter 19]{Fulton}), most of the results
on equivariant Chow groups also hold.\\

In particular, the functorial properties with respect to
proper, flat and l.c.i maps hold.\\

If $S = Spec(R)$ where $R$ is a Dedekind domain, then there
is an intersection product on $A_*^G(X)$ for $X/S$ smooth.\\

If $G$ acts freely on $X$ and a quotient $X/G$ exists,
then $A_*^G(X) = A_*(X/G)$.\\

If the stabilizer group scheme for the action of $G$ on
$X$ is finite over $S$, and a quotient $X/G$ exists,
then we expect that $A_*^G(X)_{\Q} = A_*(X/G)_{\Q}$. However,
to prove such a statement using the techniques of this paper
would require a localization long exact sequence for higher
Chow groups over an arbitrary base.\\

If $S$ is regular and $T/S$ is a split torus, then the equality of
completions of $K_i^{'T}(X)$ with respect to either the augmentation
ideal of $K_0^T(S)$ or the augmentation ideal of $K_0^T(X)$ holds (cf.
Theorem \ref{completions}). The analogous equality of completions for
(higher) Chow groups also holds.  From the torus case, we can deduce
the corresponding equality of completions for the totally split (i.e.
pulled back from the split groups over $Spec(\Z)$) classical groups $G
= Sl(n,S)$, $Sp(2n,S)$.  The argument is the same as in
Section \ref{eqrr}.  The key point is that for
these groups $K_i^{'G}(X)$ is a direct summand in $K_i^{'T}(X)$; this
can be proved (as in Lemma \ref{summand}) by realizing $G/B$ as a
sequence of projective bundles and applying Thomason's projective
bundle theorem.  (Note that for group schemes $G$, as for groups,
there is a scheme $G/B$.) For $G = SO(n,S)$, $G/B$ can be realized as
a sequence of quadric bundles, and the analysis of \cite{E-G1} applied
to deduce the result with rational
coefficients.  Once the analogue of Theorem \ref{completions} is
proved, the corresponding Riemann-Roch statements follow.\\

A form of the localization theorem for split torus
actions also holds over an arbitrary base. However, we can only
prove a localization isomorphism if the fixed locus is regularly
embedded in $X$. Again, the obstruction is the lack of
a long exact sequence for higher Chow groups over an arbitrary base.\\

Finally, if $G/S$ is smooth but not reductive and $G$ embeds
as a closed subgroup of $GL(n,S)$, then
a quotient $X \times^G U$ exists as an algebraic space. To develop
an equivariant intersection theory in this case would require further
facts about Chow groups of algebraic spaces.

\section{Appendix}
Here, we collect some useful results about actions
of algebraic groups acting on algebraic schemes in arbitrary characteristic.

\subsection{Torus actions}

\begin{lemma} \label{porb}
If $X$ is a variety with an action of a torus $T$, then there is
an open $U \subset X$ so that the stabilizer is constant for
all points of $U$.
\end{lemma}

Proof: It suffices to prove the lemma
after finite base change, so we may assume that
$T$ is split. Let $\tilde X \rightarrow X$ be the normalization
map. This map is $T$-equivariant and is an isomorphism
over an open set. Thus we may assume $X$ is normal. By Sumihiro's
theorem, the $T$ action on $X$ is locally linearizable, so it
suffices to prove the lemma when $X = V$ is a vector space and
the action is diagonal.

If $V = k^n$, then let $U = (k^*)^n$. The $n$-dimensional
torus ${\bf G}_m^n$ acts transitively on $U$ in the obvious way.
This action commutes with the given action of $T$. Thus the stabilizer
at each closed point of $U$ is the same.
\endproof

\begin{prop} \label{complement}
Let $V$ be a vector space with a linear action of a torus $T$. Let
$V^f$ be the set of points with closed orbits and trivial
stabilizers.
Then the set $V^u = V-V^f$ is a finite union of linear subspaces.
\end{prop}
Proof: Again, after finite base change we may assume that $T$ is
split.  Let $V^c$ denote the set of points in $V$ whose $T$-orbits are
closed in $V$.  We first prove that $V-V^c$ is a finite union of
linear subspaces.  Choose a basis $\{ v_i \}$ on which $T$ acts
diagonally.  If $\mbox{dim }T=1$, then the $T$-orbit of $v = \sum a_i
v_i$ is not closed if and only if the weights of the non-zero coordinates
are either all non-negative or all non-positive. Thus, $V-V^c$ is
defined by the vanishing of various subsets of coordinate hyperplanes,
hence is a finite union of linear subspaces.  For $T$ of arbitrary
dimension, $T \cdot v$ is closed iff for all 1-dimensional subtori $S
\subset T$, $S \cdot v$ is closed (this follows from \cite[Prop.
2.4]{GIT}).  This in turn holds if $S \cdot v$ is closed for a
sufficiently general $S \subset T$, so the result follows from the
case $\mbox{dim }T=1$.

To complete the proof we must show that the complement of the set
of points with trivial stabilizer is a union of linear subspaces.
This follows from two facts.\\

(1) If $G \subset T$ is a subgroup, then $L_G = \{{\bf v} \in V |
G \subset Stab({\bf v}) \}$ is a linear subspace.

(2) $V$ can be covered by a finite number of $L_G$'s by Lemma
\ref{porb}.
\endproof
\subsection{Principal bundles}
\begin{lemma} \label{q.exist} (\cite{E-G})
Let $G$ be an algebraic group. For
any $i > 0$, there is a representation $V$ of $G$ and an open set
$U \subset V$ such that $V-U$ has codimension more than $i$
and such that a principal bundle quotient $U \rightarrow U/G$
exists.
\end{lemma}
Proof. Embed $G$ into $GL(n)$ for some $n$. Assume that
$V$ is a representation
of $GL(n)$ and $U \subset V$ is an open set such that a principal
bundle quotient $U \rightarrow U/GL(n)$ exists. Since $GL(n)$
is special, this principal bundle is locally trivial in the Zariski
topology. Thus $U$ is locally isomorphic to $W \times GL(n)$ for
some open $W \subset U/GL(n)$. A quotient $U/G$ can be constructed
by patching the quotients $W \times GL(n) \rightarrow W \times (GL(n)/G)$.

We have thus reduced to the case $G=GL(n)$. Since the action of
an affine group is locally finite, there as an equivariant
closed embedding of $G \hookrightarrow V$ into a sufficiently large
vector space $V/k$. Consider the open set $U \subset V$ of points with trivial
stabilizers
which are stable for the
$G$ action on $V$. Since $G$ acts freely on itself,
$G \subset U$; hence $U$ is non-empty.
Since the stabilizers are trivial, the action on
$U$ is free, and the GIT quotient $U \rightarrow U/G$ is a principal
bundle. Now if $V_1 = V \oplus V$, then
(\cite[Proposition 1.18]{GIT}) $U_1 = (U \oplus V) \cup (V \oplus U)
\subset V_1^s$. Thus a principal bundle quotient $U_1 \rightarrow U_1/G$
exists, and the codimension of $V_1 - U_1$ is strictly smaller
than the codimension of $V-U$. Thus, by taking the direct sum of
a  sufficiently large number of copies $V$, we may assume that $V - U$
has arbitrarily high codimension. \endproof

\medskip

Let $G$ be an algebraic group, let $U$ be a scheme on which $G$ acts
freely, and suppose that a principal bundle quotient $U \rightarrow
U/G$ exists.
\begin{prop} \label{inap}
Let $X$ be an algebraic scheme with a $G$ action.
Assume that at least one of the following hypotheses holds.\\

(1) $X$ is quasi-projective with a linearized $G$-action.\\

(2) $G$ is connected and $X$ is equivariantly embedded as a closed
subscheme of a normal variety.\\

(3) $G$ is special.\\

Then a principal bundle quotient
$X \times U \rightarrow (X \times^G U)$ exists.
\end{prop}
Proof. If $X$ is quasi-projective with a linearized action, then there
is an equivariant line bundle on $X \times U$ which is relatively
ample for the projection $X \times U \rightarrow U$. By \cite[Prop 7.1]{GIT}
a principal bundle quotient $X \times^G U$ exists.

Now suppose that $X$ is normal and $G$ is connected.
By Sumhiro's theorem \cite{Sumihiro}, $X$ can be covered
by invariant quasi-projective open sets which have a linearized
$G$ action. Thus, by \cite[Prop 7.1]{GIT} we can construct
a quotient $X_G = X \times^{G} U$ by patching the quotients of
the quasi-projective open sets in the cover.

If $X$ equivariantly embeds in a normal variety $Y$, then by the above
paragraph a principal bundle quotient $Y \times U \rightarrow Y
\times^G U$ exists. Since $G$ is affine, the quotient map is affine,
and $Y \times U$ can be covered by affine invariant open sets.  Since
$X \times U$ is an invariant closed subscheme of $Y \times U$, $X
\times U$ can also be covered by invariant affines. A quotient $X
\times^G U$ can then be constructed by patching the quotients of the
invariant affines.

Finally, if $G$ is special, then $U \rightarrow U/G$ is a
locally trivial bundle in the Zariski topology. Thus
$U = \bigcup\{U_\alpha\}$ where
$\phi_{\alpha}:U_\alpha \simeq G \times W_\alpha$
for some open $W_\alpha \subset U/G$. Then $\psi_{\alpha}:
X \times U_{\alpha} \rightarrow X \times W_{\alpha}$ is a quotient,
where $\psi_{\alpha}$ is defined by the formula
$(x,w,g) \mapsto (g^{-1}x,w)$
(Here we assume that $G$ acts on the left on both factors
of $X \times G$).

\subsection{Quotients}
Following Vistoli, we define a geometric quotient $X \stackrel{\pi}
\rightarrow Y$
to be a map which satisfies properties i)-iii) of \cite[Definition
0.6]{GIT}. In particular, we do not require that ${\cal O}_Y = \pi_*({\cal O}
_X)^G$. The advantage of this definition is that is preserved under base
change. In characteristic 0 there are no inseparable extensions,
so our definition agrees with Mumford's
(\cite[Prop 0.2]{GIT}).

The following proposition is an analogue of \cite[Prop 2.6]{Vi}. The proof
is similar.
\begin{prop} \label{whizzbang}
Let $G$ act properly on a variety $X$ (hence with finite, but possibly
non-reduced stabilizers), so that a geometric quotient $X
\rightarrow Y$ exists. Then there is a commutative diagram of
quotients
$$\begin{array}{ccc} Z & \rightarrow  & X\\
\downarrow & & \downarrow\\
Q & \rightarrow & Y
\end{array}$$
where $Z \rightarrow Q$ is a principal $G$-bundle and the horizontal
maps are finite and surjective.
\end{prop}
Proof.
By \cite[Lemma p. 14]{GIT}, there is a finite map $Q \rightarrow
Y$, with $Q$ normal,
so that the pullback $X_1 \stackrel{\pi}\rightarrow Q$ has a section
in the neighborhood of every point. Cover $Q$
by a finite number of open sets $\{U_\alpha\}$ so that
$X_1 \rightarrow Q$ has a section $U_\alpha \stackrel{s_{\alpha}}
\rightarrow
V_{\alpha}$ where $V_{\alpha} = \pi^{-1}(U_{\alpha})$.

Define a $G$-map
$$\phi_{\alpha}: G \times U_\alpha \rightarrow V_\alpha$$
by the formula
$$(g,y) \mapsto gs_\alpha(y).$$
The action is proper so
each $\phi_\alpha$ is proper. Since the stabilizers
are finite, $\phi_{\alpha}$ is in fact finite.

To construct a principal bundle $Z \rightarrow Q$ we must glue
the $G \times U_{\alpha}$'s along
their intersection. To do this we will find isomorphisms
$\phi_{\alpha\beta}: s_\alpha(U_{\alpha\beta}) \rightarrow
s_{\beta}(U_{\alpha\beta})$ which satisfy the cocycle
condition.

For each $\alpha, \beta$, let $I_{\alpha\beta}$ be the scheme
which parametrizes isomorphisms of $s_\alpha$ and $s_\beta$
over $U_{\alpha\beta}$ (i.e. a section $U_{\alpha\beta}
\rightarrow I_{\alpha\beta}$ corresponds to a global isomorphism
$s_\alpha(U_{\alpha\beta}) \rightarrow s_\beta(U_{\alpha\beta})$).
The scheme $I_{\alpha\beta}$ is finite
over $U_{\alpha\beta}$ (but possibly totally ramified
in characteristic
$p$) since it is defined by the cartesian diagram
$$\begin{array}{ccc}
I_{\alpha\beta} & \rightarrow & U_{\alpha\beta} \\
\downarrow & & \small{1 \times s_\beta} \downarrow\\
G\times U_{\alpha\beta} & \stackrel{1 \times \phi_\alpha} \rightarrow &
U_{\alpha\beta} \times V_{\alpha\beta}
\end{array}$$
(Note that $I_{\alpha\alpha}$ is the stabilizer of $s_\alpha(U_\alpha)$.)

Over $U_{\alpha\beta\gamma}$ there is a composition giving
multiplication morphisms which are surjective when $\gamma = \beta$.
$$I_{\alpha\beta} \times_{U_{\alpha\beta\gamma}} I_{\beta\gamma}
\rightarrow I_{\alpha\gamma}$$ which gives
multiplication morphisms which are surjective when $\gamma = \beta$.

After a suitable finite (but possibly inseparable) base change, we may
assume that there is a section $U_{\alpha\beta} \rightarrow
I_{\alpha\beta}$ for every irreducible component of $I_{\alpha\beta}$.
(Note that $I_{\alpha\beta}$ need not be reduced.) Fix an open set
$U_{\alpha}$. For $\beta \neq \alpha$ choose a section
$\phi_{\alpha\beta}: U_{\alpha\beta} \rightarrow I_{\alpha\beta}$.
Since the $I_{\alpha\beta}$'s split completely and $I_{\alpha\alpha}$
is a group scheme, there are sections
$\phi_{\beta\alpha}:U_{\alpha\beta} \rightarrow I_{\beta\alpha}$ so
that $\phi_{\alpha\beta} \cdot \phi_{\beta\alpha}$ is the identity
section of $U_{\alpha\alpha}$. For any $\beta, \gamma$ we can define a
section of $I_{\beta\gamma}$ over $U_{\alpha\beta\gamma}$ as the
composition $\phi_{\beta\alpha} \cdot \phi_{\alpha\gamma}$. Because
$I_{\beta\gamma}$ splits, the $\phi_{\beta\alpha}$'s extend to
sections
over $U_{\beta\gamma}$.

By construction, the $\phi_{\beta\gamma}$'s satisfy the cocycle condition.
We can now define $Z$ by gluing the sets $G \times U_{\beta}$
along the $\phi_{\beta\gamma}$'s.
\endproof


\begin{thebibliography}{99}
\bibitem[At]{Atiyah} M. Atiyah, {\it Bott periodicity and the index of
elliptic operators}, Quart. J. Math. Oxford (2) {\bf 19} (1968),
113-40.
\bibitem[A-B]{A-B} M. Atiyah, R. Bott, {\it
The moment map and equivariant cohomology}, Topology {\bf 23} (1984),
1-28.
\bibitem[A-S]{A-S} M. Atiyah, G. Segal, {\it Equivariant K-theory and
completion }, J. Diff. Geom. {\bf 3} (1969), 1-18.
\bibitem[B-V]{B-V} N. Berline, M. Vergne, {\it Classes caract\'eristiques
equivariantes. Formule de localization en cohomologie \'equivariante},
C.R. Acad. Sc. Paris {\bf 295} (1982), 539-541.
\bibitem[Bl]{Bl} S. Bloch, {\it Algebraic cycles and higher
$K$-theory}, Adv. Math. {\bf 61} (1986), 267-304, and
{\it The moving lemma for higher Chow groups}, J. Alg. Geom. {\bf 3} (1994),
537-568.
\bibitem[Bo]{Borel} A. Borel,
{\it Linear Algebraic Groups}, 2nd enlarged edition,
Graduate Texts in Mathematics {\bf 126}, Springer Verlag (1991).
\bibitem[Br]{Br} R.E. Briney, {\it Intersection theory on
quotients of algebraic varieties}, Am. J. Math. {\bf 84} (1962),
217-238.
\bibitem[CEPT]{CEPT} T. Chinburg, B. Erez, G. Pappas, M.J. Taylor,
{\it Riemann-Roch type theorems for arithmetic schemes with a
finite group action}, preprint.
\bibitem[D-M]{D-M} P. Deligne, D. Mumford, {\it Irreducibility of
the space of curves of a given genus}, Publications I.H.E.S.
{\bf 36} (1969), 75-109.
\bibitem[SGA3]{SGA3} M. Demazure, A. Grothendieck, {\it Sch\'emas en
groupes}, Springer Lecture Notes in Math, {\bf 151-153} (1970).
\bibitem[Ed]{Edidin} D. Edidin, {\it The codimension-two homology of the moduli
space of
stable curves is algebraic}, Duke Math. J. {\bf 67} (1992), 241-272.
\bibitem[E-G]{E-G} D. Edidin, W. Graham, {\it Characteristic
classes in the Chow ring}, preprint.
\bibitem[E-G1]{E-G1} D. Edidin, W. Graham, {\it Characteristic classes
and quadric bundles}, Duke Math. Journal, {\bf 78} (1995), 277-299.
\bibitem[E-S]{E-S} G. Ellingsrud, S. Stromme, {\it Bott's formula and
enumerative geometry}, J. Amer. Math. Soc. {\bf 9} (1996), 175-194.
\bibitem[Fu]{Fulton} W. Fulton, {\it Intersection Theory},
Ergebnisse, 3. Folge, Band 2, Springer Verlag (1984).
\bibitem[F-L]{F-L} W. Fulton, S. Lang, {\it Riemann-Roch algebra},
Springer-Verlag (1985).
\bibitem[Gr]{Grayson} D. Grayson, {\it Products in $K$-theory and intersecting
cycles}, Inv. Math. {\bf 47} (1978), 71-83.
\bibitem[Gi]{Gi} H. Gillet, {\it Intersection theory on algebraic
stacks and $Q$-varieties}, J. Pure Appl. Alg., {\bf 34} (1984),
193-240.
\bibitem[GIT]{GIT} D. Mumford, J. Fogarty, F. Kirwan, {\it Geometric
Invariant Theory}, 3rd enlarged edition, Springer-Verlag (1994).
\bibitem[H-Y]{H-Y} R. H\"ubl, A. Yekutieli, {\it Adelic Chern
forms and the Bott residue formula}, preprint.
\bibitem[Iv]{Iv} B. Iversen, {\it A fixed point formula for actions
of tori on algebraic varieties}, Inv. Math. {\bf 16} (1972), 229-236.
\bibitem[I-N]{I-N} B. Iversen, H. Nielsen, {\it Chern numbers and
diagonalizable groups}, J. London Math. Soc. {\bf 11} (1975), 223-232.
\bibitem[Ki]{Kimura} S. Kimura, {\it Fractional intersection
and bivariant theory}, Comm. Alg. {\bf 20} (1992) 285-302.
\bibitem[Ko]{kock} B. K\"ock, {\it The Grothendieck Riemann-Roch theorem
in the higher $K$-theory of group scheme actions}, preprint.
\bibitem[Mu]{Mu} D. Mumford, {\it Towards an enumerative geometry of
the moduli space of curves}, Prog. Math. {\bf 36}, (1983), 271-328.
\bibitem[Sem-Chev]{Sem-Chev} {\it Anneau de Chow et applications}, Seminaire
Chevalley, Secr\'etariat math\'ematique, Paris (1958).
\bibitem[Se]{Seshadri} C.S. Seshadri, {\it Geometric reductivity over
arbitrary base}, Adv. Math {\bf 26} (1977), 225-274.
\bibitem[Su]{Sumihiro} H. Sumihiro, {\it Equivariant completion II},
J. Math. Kyoto {\bf 15} (1975), 573-605.
\bibitem[Th]{Thomason} R. Thomason, {\it Algebraic $K$-theory
of group scheme actions}, in {\it Algebraic topology and
algebraic $K$-theory} (W. Browder, editor), Annals of Math Studies
{\bf 113} (1987), 539-563.
\bibitem[Th1]{Th1} R. Thomason, {\it Comparison of equivariant
algebraic and topological $K$-theory}, Duke Math. J.
{\bf 53} (1986), 795-825.
\bibitem[Th2]{Th2} R. Thomason, {\it Lefschetz-Riemann-Roch theorem
and coherent trace formula}, Inv. Math. {\bf 85} (1986), 515-543.
\bibitem[To]{To} B. Totaro, {\it The Chow ring of the symmetric
group}, preprint.
\bibitem[Vi]{Vi} A. Vistoli, {\it Intersection theory on algebraic
stacks and their moduli}, Inv. Math. {\bf 97} (1989), 613-670.
\bibitem[Vi1]{Vi3} A. Vistoli, {\it Equivariant Grothendieck groups and
equivariant Chow groups} in {\it Classification of
irregular varieties (Trento, 1990)}, Lecture Notes in Math.
{\bf 1515}, 112--133, Springer, Berlin, 1992.
\end{thebibliography}
\end{document}